%% file: ms.tex
\newcount\Comments  
\RequirePackage{fix-cm}
\documentclass[final,leqno,onefignum,onetabnum]{siamltex}

\usepackage{graphicx}
\usepackage[utf8]{inputenc}
\usepackage[english]{babel}
\usepackage[automark]{scrpage2}
\usepackage{amsfonts}
\usepackage{amsmath}
\usepackage{amssymb}
\usepackage{color}
\usepackage{graphicx}
\usepackage{sidecap}
\usepackage{multirow}
\usepackage{booktabs}
\usepackage{bm}

\usepackage{tabularx}
\usepackage[font=small]{caption}
\usepackage{subfig}

\usepackage{longtable}
\usepackage{rotating}

\usepackage{algorithm}
\usepackage{algpseudocode}

\graphicspath{{figures/}}
\DeclareGraphicsExtensions{.pdf,.eps,.png,.jpg,.jpeg}

\newcommand{\Zcal}{\mathcal{Z}} 
\newcommand{\Ycal}{\mathcal{Y}} 

\newcommand{\bfc}{\boldsymbol c} 
\newcommand{\bfe}{\boldsymbol e} 
\newcommand{\bfr}{\boldsymbol r} 
\newcommand{\bfy}{\boldsymbol y} 
\newcommand{\bfz}{\boldsymbol z} 

\newcommand{\bfSigma}{\boldsymbol \Sigma} 
\newcommand{\bfeps}{\boldsymbol \epsilon} 

\newcommand{\hi}{\text{hi}}
\newcommand{\lo}{\text{lo}}
\newcommand{\fhigh}{f_{\hi}}
\newcommand{\fl}{f_{\lo}}
\newcommand{\chigh}{c_{\hi}}
\newcommand{\clow}{c_{\lo}}

\newcommand{\dvComp}{z} 
\newcommand{\dv}{\boldsymbol \dvComp} 
\newcommand{\ov}{\bfy} 
\newcommand{\nrpts}{m} 
\newcommand{\obs}{\ov_{\text{obs}}} 
\newcommand{\Var}{\mathbb{V}}
\newcommand{\MCEst}{\bar{s}}

\newcommand{\kibitz}[2]{\ifnum\Comments=1\textcolor{#1}{#2}\fi}

\ihead[]{}
\ohead[]{}
\ifoot[]{}
\ofoot[]{}

\title{Survey of multifidelity methods in uncertainty propagation, inference, and optimization}

\author{Benjamin Peherstorfer\thanks{Department of Mechanical Engineering and Wisconsin Institute for Discovery, University of Wisconsin-Madison, Madison, WI 53706} \and Karen Willcox\thanks{Department of Aeronautics \& Astronautics, MIT, Cambridge, MA 02139} \and Max Gunzburger\thanks{Department of Scientific Computing, Florida State University, Tallahassee, FL 32306-4120}}

\begin{document}
\maketitle
\slugger{mms}{xxxx}{xx}{x}{x--x}

\begin{abstract}
In many situations across computational science and engineering, multiple computational models are available that describe a system of interest. These different models have varying evaluation costs and varying fidelities. Typically, a computationally expensive high-fidelity model describes the system with the accuracy required by the current application at hand, while lower-fidelity models are less accurate but computationally cheaper than the high-fidelity model. Outer-loop applications, such as optimization, inference, and uncertainty quantification, require multiple model evaluations at many different inputs, which often leads to computational demands that exceed available resources if only the high-fidelity model is used. This work surveys multifidelity methods that accelerate the solution of outer-loop applications by combining high-fidelity and low-fidelity model evaluations, where the low-fidelity evaluations arise from an explicit low-fidelity model (e.g., a simplified physics approximation, a reduced model, a data-fit surrogate, etc.) that approximates the same output quantity as the high-fidelity model. The overall premise of these multifidelity methods is that low-fidelity models are leveraged for speedup while the high-fidelity model is kept in the loop to establish accuracy and/or convergence guarantees. We categorize multifidelity methods according to three classes of strategies: adaptation, fusion, and filtering. The paper reviews multifidelity methods in the outer-loop contexts of uncertainty propagation, inference, and optimization.
\end{abstract}



\pagestyle{myheadings}
\thispagestyle{plain}
\markboth{PEHERSTORFER, WILLCOX, AND GUNZBURGER}{MULTIFIDELITY METHODS}

\input{intro}

\input{mm}

\input{uq}

\input{invprob}

\input{opti}

\input{outlook}

\section*{Acknowledgements}
The first two authors acknowledge support of the AFOSR MURI on multi-information sources of multi-physics systems under Award Number FA9550-15-1-0038, the United States Department of Energy Applied Mathematics Program, Awards DE-FG02-08ER2585 and \text{DE-SC0009297}, as part of the DiaMonD Multifaceted Mathematics Integrated Capability Center, DARPA EQUiPS Award UTA15-001067, and the MIT-SUTD International Design Center. The third author was supported by the US Department of Energy Office of Science grant \text{DE-SC0009324} and US Air Force Office of Research grant FA9550-15-1-0001.

\bibliography{multi}
\bibliographystyle{abbrv}

\end{document}

%% file: intro.tex
\section{Introduction}
\label{sec:Intro}
We begin by introducing the setting and concepts surveyed in this paper:
Section~\ref{sec:intro-mf-models} defines the setting of multifidelity models and Section~\ref{sec:intro-mf} introduces the concepts of multifidelity methods. Section~\ref{sec:intro:models} discusses different types of low-fidelity models that may arise in the multifidelity setting. Section~\ref{sec:intro-ol} defines the three outer-loop applications of interest: uncertainty propagation, statistical inference, and optimization. Section~\ref{sec:intro-outline} outlines the remainder of the paper.

\subsection{Multifidelity models} \label{sec:intro-mf-models}
Models serve to support many aspects of computational science and engineering, from discovery to design to decision-making and more. In some of these settings, one primary purpose of a model is to characterize the input-output relationship of the system of interest---the input describes the relevant system properties and environmental conditions, and the output describes quantities of interest to the task at hand.
In this context, evaluating a model means performing a numerical simulation that implements the model, computes a solution, and thus maps an input onto an approximation of the output. For example, the numerical simulation might involve solving a partial differential equation (PDE), or solving a system of ordinary differential equations, or applying a particle method. Mathematically, we denote a model as a function $f: \Zcal \to \Ycal$ that maps an input $\dv \in \Zcal$ to an output $\ov \in \Ycal$, where $\Zcal \subseteq \mathbb{R}^d$ is the domain of the inputs of the model, with dimension $d \in \mathbb{N}$, and $\Ycal \subseteq \mathbb{R}^{d^{\prime}}$ is the domain of the outputs of the model, with dimension $d^{\prime} \in \mathbb{N}$.
Model evaluations (i.e., evaluations of $f$) incur computational costs $c \in \mathbb{R}_+$ that typically increase with the accuracy of the approximation of the output, where $\mathbb{R}_+ = \{x \in \mathbb{R} \, : \, x > 0\}$ is the set of positive real numbers.

In many situations, multiple models are available that estimate the same output quantity with varying approximation qualities and varying computational costs. We define a {\em high-fidelity model} $\fhigh: \Zcal \to \Ycal$ as a model that estimates the output with the accuracy that is necessary for the current task at hand. We define a {\em low-fidelity model} $\fl: \Zcal \to \Ycal$ as a model that estimates the same output with a lower accuracy than the high-fidelity model.
The costs $\chigh \in \mathbb{R}_+$ of the high-fidelity model $\fhigh$ are typically higher than the costs $\clow \in \mathbb{R}_+$ of a low-fidelity model $\fl$. More generally, we consider $k \in \mathbb{N}$ low-fidelity models, $\fl^{(1)}, \dots, \fl^{(k)}$, that each represent the relationship between the input and the output, $\fl^{(i)}: \Zcal \to \Ycal, \ i=1,\dots,k$, and we denote the cost of evaluating model $\fl^{(i)}$ as $\clow^{(i)}$.

\subsection{Multifidelity methods for the outer loop} \label{sec:intro-mf}

The use of principled approximations to accelerate computational tasks has long been a mainstay of scalable numerical algorithms. For example,  quasi-Newton optimization methods \cite{DavidonQNewton,FletcherPaulQNewton,broyden_class_1965} construct approximations of Hessians and apply low-rank updates to these approximations during the Newton iterations. Solvers based on Krylov subspace methods \cite{lanczos_iteration_1950,arnoldi_principle_1951,lanczos_solution_1952,saad_gmres:_1986} and on Anderson relaxation \cite{Anderson:1965:IPN:321296.321305,WalkerAnderson,AndersonConv} perform intermediate computations in low-dimensional subspaces that are updated as the computation proceeds. Whereas these methods---and many others across the broad field of numerical algorithms---embed principled approximations within a numerical solver, we focus in this paper on the particular class of multifidelity methods that invoke explicit approximate models in solution of an outer-loop problem. We define this class of methods more precisely below; first we introduce the notion of an outer-loop application problem.

We use the term {\em outer-loop application} to define computational applications that form outer loops around a model---where in each iteration an input $\dv \in \Zcal$ is received and the corresponding model output $f(\dv)$ is computed, and an overall outer-loop result is obtained at the termination of the outer loop. For example, in optimization, the optimizer provides at each iteration the design variables to evaluate (the input) and the model must evaluate the corresponding objective function value, the constraints values, and possibly gradient information (the outputs). At the termination, an optimal design is obtained (the outer-loop result). Another outer-loop application is uncertainty propagation, which can be thought of conceptually as a loop over realizations of the input, requiring the corresponding model evaluation for each realization. In uncertainty propagation, the outer-loop result is the estimate of the statistics of interest. Other examples of outer-loop applications include inverse problems, data assimilation, control problems, and sensitivity analysis.\footnote{We adopted the term ``outer loop," which is used by a number of people in the community, although it does not appear to have been formally defined.  Ref.~\cite{Keyes2011} gives specific examples of outer-loop applications in the context of petroleum production. Ref.~\cite[Chapter~10.1]{ExaReport} discusses outer-loop applications in uncertainty quantification and optimization. } Note that although it is helpful for the exposition to think of outer-loop applications as loops, they are often \emph{not} implemented as such. For example, in uncertainty propagation, once the realizations of the input have been drawn, the model outputs can be typically computed in parallel.

The term {\em many-query application} is often used to denote applications that evaluate a model many times \cite{RozzaPateraSurvey}, a categorization that applies to most (if not all) outer-loop applications. We distinguish between many-query and outer-loop applications by considering the latter to be the class of applications that target a specific outer-loop result. In contrast, many-query applications do not necessarily target a specific outer-loop result (and thus the set of outer-loop applications is essentially a subset of the set of many-query applications). For example, performing a parameter study is many-query but does not necessarily lead to a specific outer-loop result. This distinction is important in the discussion of multifidelity methods, since accuracy and/or convergence will be assessed relative to a specific outer-loop result.

The accuracy of the outer-loop result, as required by the problem at hand, can be achieved by using the high-fidelity model $\fhigh$ in each iteration of the outer loop; however, evaluating the high-fidelity model in each iteration often leads to computational demands that exceed available resources. Simply replacing the high-fidelity model $\fhigh$ with a low-fidelity model $\fl$ can result in significant speedups but leads to a lower---and typically unknown---approximation quality of the outer-loop result. This is clearly unsatisfactory and motivates the need for multifidelity methods.

\begin{figure}
\centering
\subfloat[][single-fidelity approach with\\ high-fidelity model]{\includegraphics[width=0.33\columnwidth]{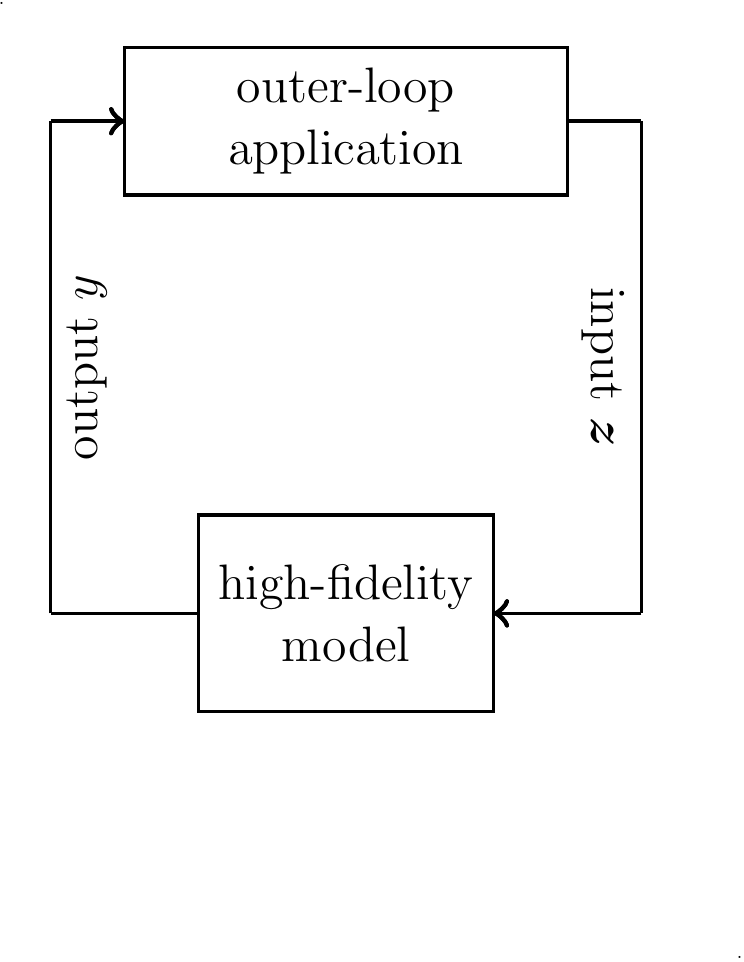}}
\subfloat[][single-fidelity approach with\\ low-fidelity model]{\includegraphics[width=0.33\columnwidth]{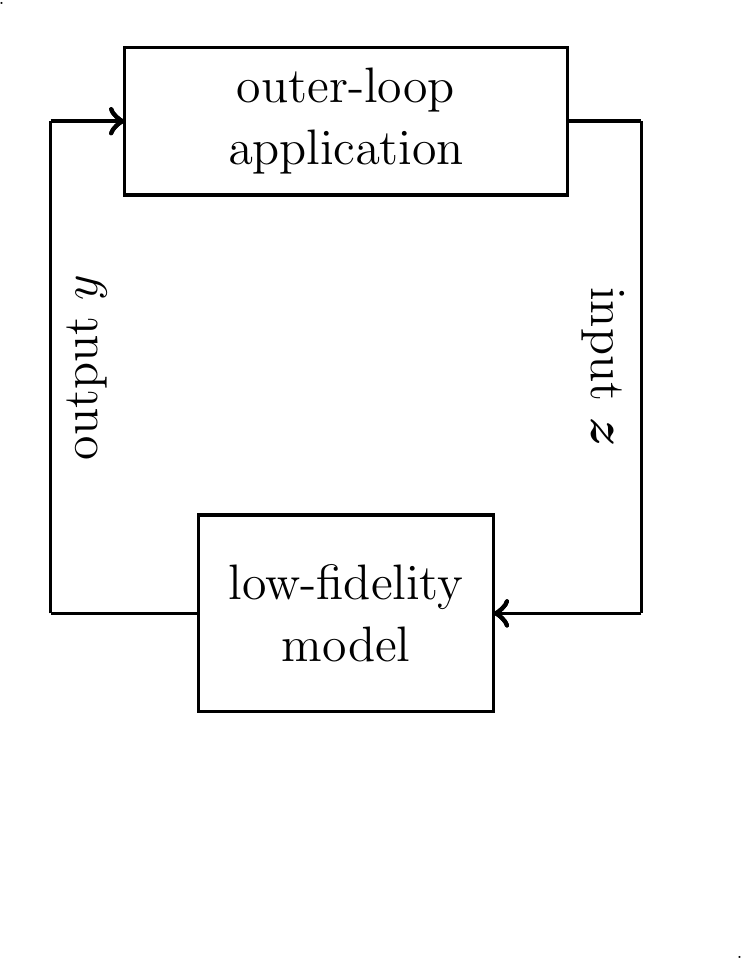}}
\subfloat[][multifidelity approach with high-fidelity model and multiple low-fidelity models]{\includegraphics[width=0.33\columnwidth]{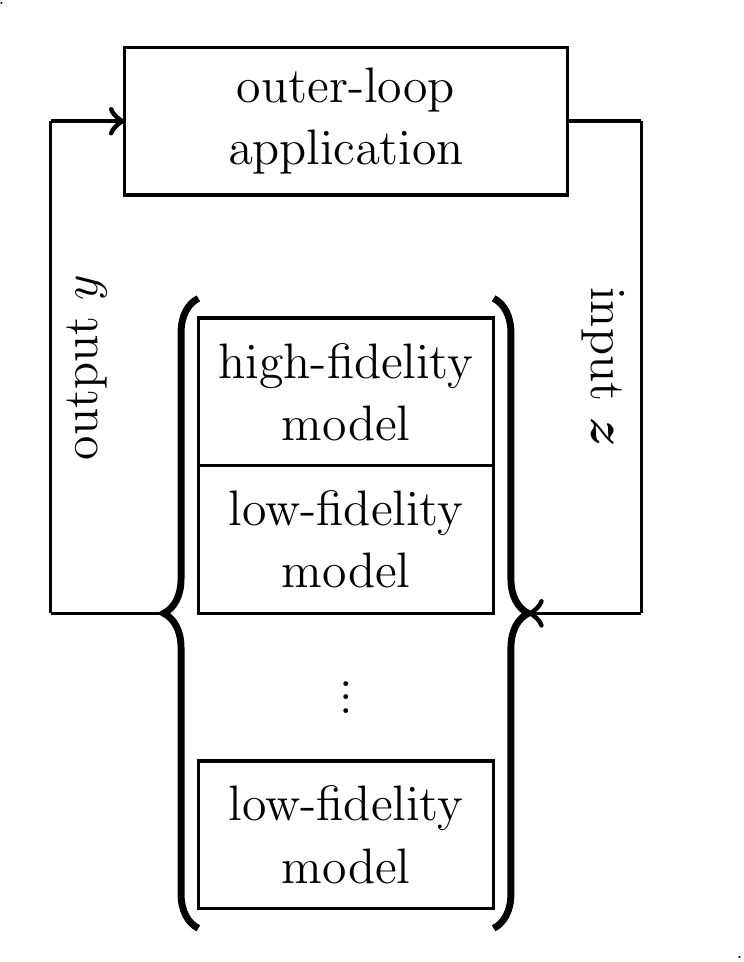}}
\caption{Multifidelity methods combine the high-fidelity model with low-fidelity models. The low-fidelity models are leveraged for speedup and the high-fidelity model is kept in the loop to establish accuracy and/or convergence guarantees on the outer-loop result.}
\label{fig:Multifidelity}
\end{figure}

We survey here multifidelity methods for outer-loop applications. We consider the class of multifidelity methods that have two key properties: (1)~They leverage a low-fidelity model $\fl$ (or in the general case multiple low-fidelity models $\fl^{(1)}, \dots, \fl^{(k)}, k \in \mathbb{N}$), to obtain computational speedups, and (2)~they use recourse to the high-fidelity model $\fhigh$ to establish accuracy and/or convergence guarantees on the outer-loop result, see Figure~\ref{fig:Multifidelity}.
Thus, \emph{multifidelity methods} use low-fidelity models to reduce the runtime where possible, but recourse to the high-fidelity model to preserve the accuracy of the outer-loop result that would be obtained with a method that uses only the high-fidelity model. The two key ingredients of multifidelity methods are (1)~low-fidelity models $\fl^{(1)}, \dots, \fl^{(k)}$, that provide useful approximations of the high-fidelity model $\fhigh$, and (2)~a \emph{model management} strategy that distributes work among the models while providing theoretical guarantees that establish the accuracy and/or convergence of the outer-loop result.

Note that a crucial component of this characterization of multifidelity methods for outer-loop problems is the use of explicit low-fidelity models that approximate the same output quantity as the high-fidelity model. This distinguishes the methods from those that embed approximations within the solver itself, such as the quasi-Newton and Krylov subspace methods discussed above.

\begin{figure}
\centering
\includegraphics[width=0.5\columnwidth]{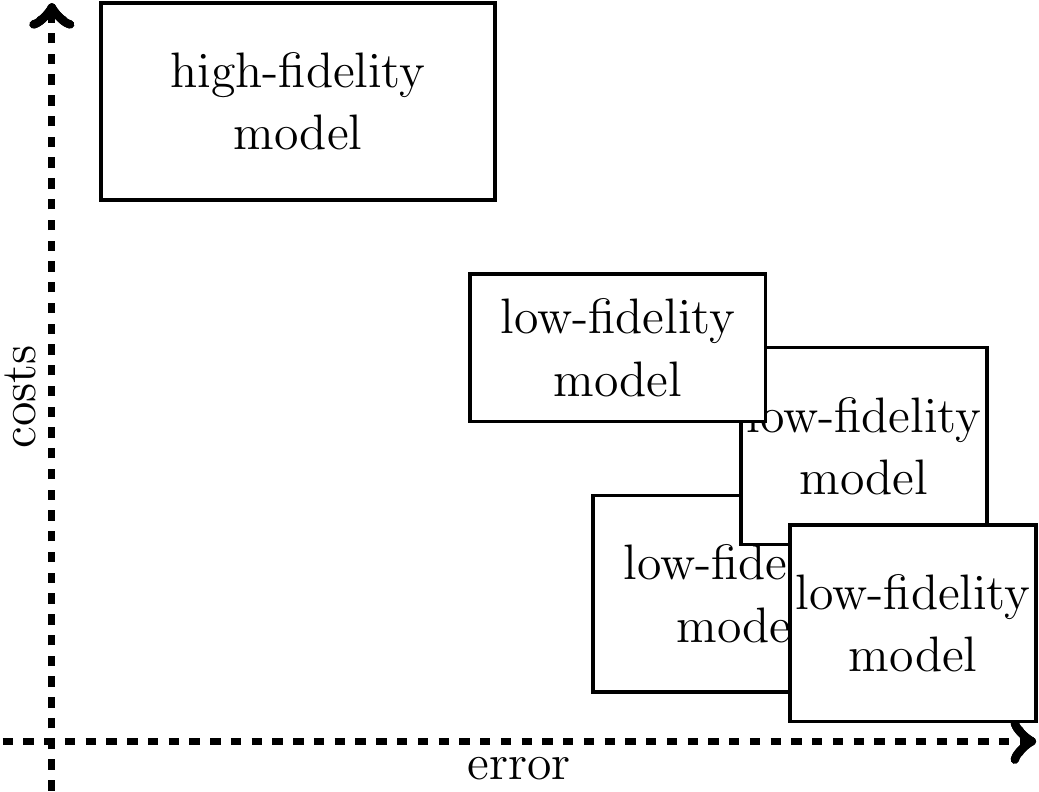}
\caption{In many situations, different types of low-fidelity models are available, e.g., coarse-grid approximations, projection-based reduced models, data-fit interpolation and regression models, machine-learning-based models, and simplified models. The low-fidelity models vary with respect to error and costs. Multifidelity methods leverage these heterogeneous types of low-fidelity models for speedup.}
\label{fig:Models}
\end{figure}

The multifidelity methods we survey are applicable to a broad range of problems, but of particular interest is the setting of a high-fidelity model that corresponds to a fine-grid discretization of a PDE that governs the system of interest. In this setting, coarse-grid approximations have long been used as cheaper approximations. Varying the discretization parameters generates a hierarchy of low-fidelity models. We are here interested in richer and more heterogeneous sets of models, including projection-based reduced models \cite{SirovichMethodOfSnapshots,RozzaPateraSurvey,gugercin_2008,SIREV}, data-fit interpolation and regression models \cite{ForresterKriging,forrester_recent_2009}, machine-learning-based models such as support vector machines (SVMs) \cite{VapnikBook,SVM,CC01a}, and other simplified models \cite{majda_quantifying_2010,ng_monte_2015},
see Figure~\ref{fig:Models}. We further discuss types of low-fidelity models in Section~\ref{sec:intro:models}. In a broader sense, we can think of the models as information sources that describe the input-output relationships of the system of interest. In that broader sense, expert opinions, experimental data, and historical data are potential information sources. We restrict the following discussion to models, because all of the multifidelity methods that we survey are developed in the context of models; however, we note that many of these multifidelity methods could potentially be extended to this broader class of information sources. 

\begin{figure}
\centering
\includegraphics[width=1.0\columnwidth]{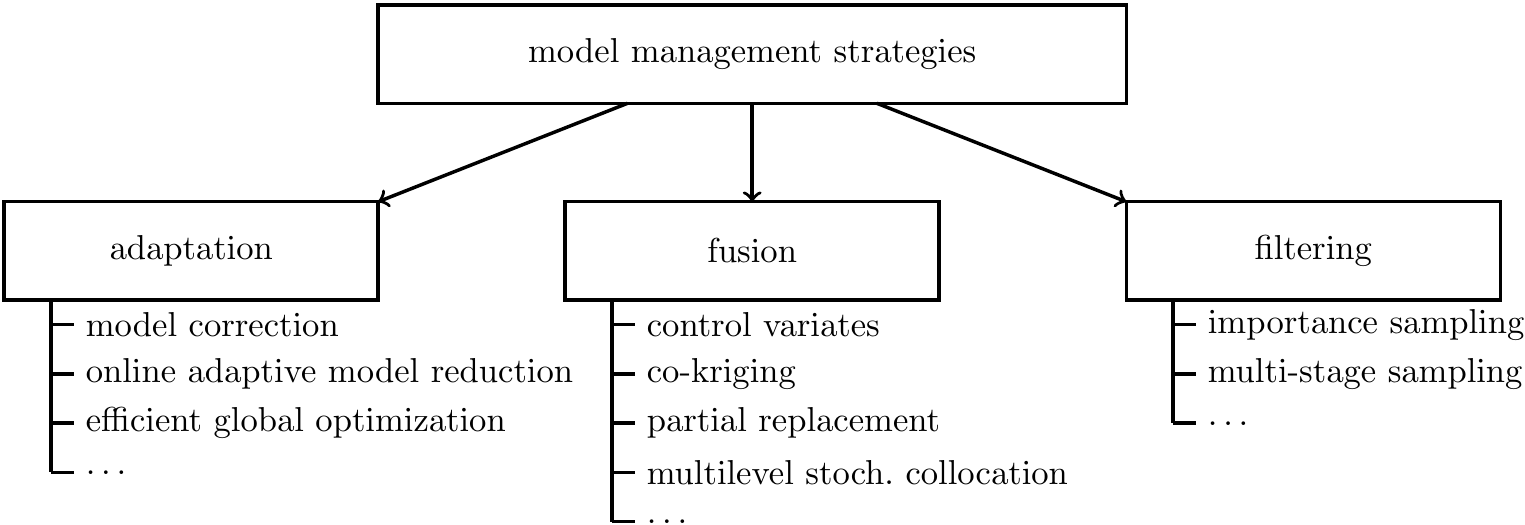}
\caption{We distinguish between three model management strategies: adaptation, fusion, and filtering.}
\label{fig:MMTypes}
\end{figure}

Model management serves two purposes. First is to balance model evaluations among the models (i.e., to decide which model to evaluate when). Second is to guarantee the same accuracy in the outer-loop result as if only the high-fidelity model were used.
We distinguish between three types of model management strategies (see Figure~\ref{fig:MMTypes}): (1)~{\em adapting} the low-fidelity model with information from the high-fidelity model,   (2)~{\em fusing} low- and high-fidelity model outputs, and (3)~{\em filtering} to use the high-fidelity model only when indicated by a low-fidelity filter.\footnote{Note that we use the term {\em filter} to denote selective evaluation based on the low-fidelity model. This differs from the predominant usage in signal processing and uncertainty quantification, where filtering describes the estimation of the state of a dynamical system from noisy and incomplete data (e.g., Kalman filter, particle filter).} The appropriate model management strategy for the task at hand typically depends on the nature of the outer-loop application. We survey model management techniques that fall into these three categories in Section~\ref{sec:MM:strategies}.

\paragraph{Comparison to multilevel methods}
Multilevel methods have a long history in computational science and engineering, e.g., multigrid methods \cite{brandt_multi-level_1977,HackbuschMultigrid,MultigridTutorial,McCormick1987,trottenberg01}, multilevel preconditioners \cite{bramble_parallel_1990,DahmenPreconditioner}, and multilevel function representations \cite{yserentant_hierarchical_1986,Bank1988,Wavelets,bungartz_sparse_2004}. Multilevel methods typically derive a hierarchy of low-fidelity models of the high-fidelity model by varying a parameter. For example, the parameter could be the mesh width and thus the hierarchy of low-fidelity models would be the hierarchy of coarse-grid approximations. A common approach in multilevel methods is to describe the approximation quality and the costs of the low-fidelity model hierarchy with rates, and then to use these rates to distribute work among the models. In this paper, we consider more general low-fidelity models with properties that cannot necessarily be well described by rates. Even though many multilevel methods are applicable to more heterogeneous models than coarse-grid approximations, describing the model properties by rates only, and consequently distributing work with respect to rates, can be too coarse a description and can miss important aspects of the models.  Furthermore, in our setting, low-fidelity models are often given and cannot be easily generated on request by varying a (e.g., discretization) parameter. The multifidelity techniques that we describe here explicitly take such richer sets of models into account.

\paragraph{Comparison to traditional model reduction}
Traditionally, model reduction \cite{AntoulasBook,RozzaPateraSurvey,SIREV} first constructs a low-fidelity reduced model, and then replaces the high-fidelity model with the reduced model in an outer-loop application. Replacing the high-fidelity model often leads to significant speedups, but it also means that the accuracy of the outer-loop result depends on the accuracy of the reduced model. In some settings, error bounds or error estimates are available for the reduced model outputs \cite{RozzaPateraSurvey,FLD:FLD867,grepl_efficient_2007}, and it may be possible to translate these error estimates on the model outputs into error estimates on the outer loop result. In contrast, multifidelity methods establish accuracy and convergence guarantees---instead of providing error bounds and error estimates only---by keeping the high-fidelity model in the loop and thus trading some speedup for guarantees---even if the quality of the low-fidelity model is unknown.

\begin{figure}
\centering
\includegraphics[width=1.0\columnwidth]{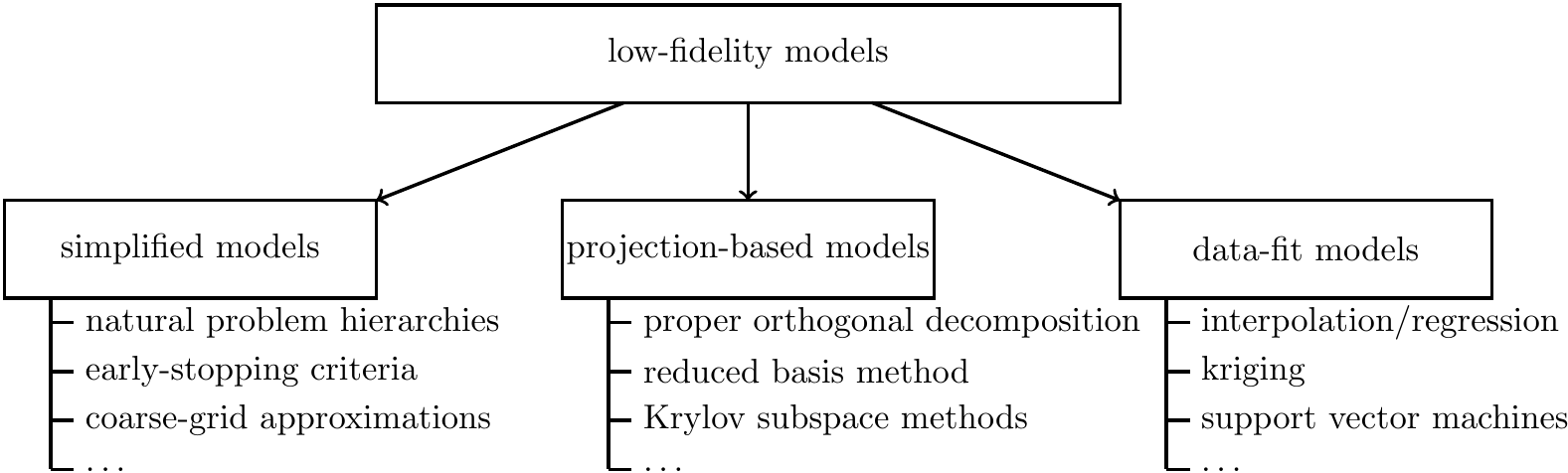}
\caption{We categorize low-fidelity models as being of three types: simplified models, projection-based models, and data-fit models.}
\label{fig:ModelTypes}
\end{figure}

\subsection{Types of low-fidelity models} \label{sec:intro:models}
We categorize low-fidelity as being of three types: simplified low-fidelity models, projection-based low-fidelity models, and data-fit low-fidelity models. Figure~\ref{fig:ModelTypes} depicts this categorization. For a given application, knowledge of and access to the high-fidelity model affects what kind of low-fidelity models can be created. In some cases, the high-fidelity system has a known structure that can be exploited to create low-fidelity models. In other cases, the high-fidelity models are considered to be ``black box": they can be evaluated at the inputs in $\Zcal$ to obtain outputs in $\Ycal$, however no details are available on how the outputs are computed.

\paragraph{Simplified low-fidelity models} Simplified models are derived from the high-fidelity model by taking advantage of domain expertise and in-depth knowledge of the implementation details of the high-fidelity model. Domain expertise allows the derivation of several models with different computational costs and fidelities that all aim to estimate the same output of interest of the system. For example, in computational fluid dynamics, there is a clear hierarchy of models for analyzing turbulent flow. From high to low fidelity, these are direct numerical simulations (DNS), large eddy simulations (LES), and Reynolds averaged Navier-Stokes (RANS). All these model turbulent flows, but DNS resolves the whole spatial and time domain to the scale of the turbulence, LES eliminates small scale behavior, and RANS applies the Reynolds decomposition to average over time. In aerodynamic design, an often-employed hierarchy comprises RANS, Euler equations, and potential theory \cite{Han2013177}. The supersonic aerodynamic design problem in \cite{choi_multifidelity_2008,Choi2009} employs the Euler equations, a vortice lattice model, and a classical empirical model. Similar hierarchies of models exist in other fields of engineering. Models for subsurface flows through karst aquifers reach from simple continuum pipe flow models \cite{NUM:NUM20579} to coupled Stokes and Darcy systems \cite{CoupledStokesDarcy}. In climate modeling, low-fidelity models consider only a limited number of atmospheric effects whereas high-fidelity models are fully-coupled atmospheric and oceanic simulation models \cite{held_gap_2005,majda_quantifying_2010}.
There are more general concepts to derive low-fidelity models by simplification, which also require domain expertise but which are applicable across disciplines. Coarse-grid discretizations are an important class of such approximations.
As another example, in many settings a low-fidelity model can be derived by neglecting nonlinear terms. For example, lower-fidelity linearized models are common in aerodynamic and structural analyses \cite{piperni_development_2013}. Yet another example is when the high-fidelity model relies on an iterative solver (e.g., Krylov subspace solvers or Newton's method), a low-fidelity model can be derived by loosening the residual tolerances of the iterative method---thus, to derive a low-fidelity approximation, the iterative solver is stopped earlier than if a high-fidelity output were computed.

\paragraph{Projection-based low-fidelity models} Model reduction derives low-fidelity models from a high-fidelity model by mathematically exploiting the problem structure, rather than using domain knowledge of the problem at hand. These methods proceed by identifying a low-dimensional subspace that is constructed so as to retain the essential character of the system input-output map. Projecting the governing equations onto the low-dimensional subspace yields the reduced model. The projection step is generally (but not always) intrusive and requires knowledge of the high-fidelity model structure. There are a variety of ways to construct the low-dimensional subspace, see \cite{SIREV} for a detailed review. One common method is the proper orthogonal decomposition (POD) \cite{SirovichMethodOfSnapshots,PODFluid,rathinam_new_2003,KunischPOD,Kunisch2001}, which uses so-called ``snapshots"---state vectors of the high-fidelity model at selected inputs---to construct a basis for the low-dimensional subspace.
POD is a popular basis generation method because it is applicable to a wide range of problems, including time-dependent and nonlinear problems \cite{deim2010,PODGreedy,Kunisch2001,KunischPOD,KutzLibrary}. Another basis generation approach is based on centroidal Voronoi tessellation (CVT) \cite{MaxCVTReview}, where a special Voronoi clustering of the snapshots is constructed. The reduced basis is then derived from the generators of the Voronoi clustering. The work \cite{Burkardt2006337} discusses details on CVT-based basis construction. A combination of POD and CVT-based basis construction is introduced in \cite{PODCVTcombinationDu2003}. There are also methods based on Krylov subspaces to generate a reduced basis \cite{FeldmannKrylov,gallivan}, including multivariate Pad\'{e} approximations and tangential interpolation for linear systems \cite{bai_krylov_2002,BaurInterp,freund_model_2003,gugercin_2008}.
Dynamic mode decomposition is another basis generation method that is popular in the context of computational fluid dynamics \cite{SchmidDMD,Tu2014391,DMDWithControl}. Balanced truncation \cite{moore_principal_1981,mullis_synthesis_1976} is a common basis construction method used in the systems and control theory community. For stable linear time-invariant systems, balanced truncation provides a basis that guarantees asymptotically stable reduced systems and provides an error bound \cite{AntoulasBook,SerkanSurvey}. Another basis generation approach is the reduced basis method \cite{RozzaPateraSurvey,grepl_efficient_2007,Rozza20071244,Rozza2013}, where orthogonalized carefully selected snapshots are the basis vectors.
Depending on the problem of interest, these reduced basis models can be equipped with cheap \emph{a posteriori} error estimators for the reduced model outputs \cite{HUYNH2007473,grepl_efficient_2007,refId0,RozzaPateraSurvey,Urban2012203,MasaRBM}. Efficient error estimators can also sometimes be provided for other basis generation methods, such as the POD \cite{Hinze2004,HinzeVolkwein2008}.

\paragraph{Data-fit low-fidelity models} Data-fit low-fidelity models are derived directly from inputs and the corresponding outputs of the high-fidelity model. Thus, data-fit models can be derived from black-box high-fidelity models because only the inputs and outputs of the high-fidelity model need to be available. In many cases, data-fit models are represented as linear combinations of basis functions. Data-fit models are constructed by fitting the coefficients of the linear combination via interpolation or regression to the inputs and the corresponding high-fidelity model outputs. The choice of the interpolation and regression bases is critical for the approximation quality of the data-fit models. Polynomials, e.g., Lagrange polynomials, are classical basis functions that can be used for deriving data-fit models. Piecewise-polynomial interpolation approaches allow use of low-degree polynomials, which avoids problems with global polynomial interpolation of high degree, e.g., Runge's phenomenon. If the inputs are low dimensional, a multi-variate data-fit model can be derived with tensor product approaches. In higher-dimensional settings, discretization methods based on sparse grids \cite{bungartz_sparse_2004} can be employed. Radial basis functions are another type of basis functions that are widely used for constructing data-fit models \cite{RasmussenKriging,forrester_recent_2009}. If based on the Gaussian density function, radial basis functions typically lead to numerically stable computations of coefficients of the linear combination of the data-fit model. Often the radial basis functions depend on hyper-parameters, e.g., the bandwidth of the Gaussian density function. Well-chosen hyper-parameters can greatly improve the approximation accuracy of the data-fit model but the optimization for these hyper-parameters is often computationally expensive \cite{forrester_recent_2009}. A widely used approach to interpolation with radial basis functions is kriging, for which a sound theoretical understanding has been obtained and efficient approaches to optimize for the hyper-parameters have been developed \cite{matheron_principles_1963,sacks_design_1989,keane_wing_2003,martin_use_2005,RasmussenKriging,ForresterKriging}. 
In particular, kriging models are equipped with error indicators, see, e.g., \cite{sacks_design_1989}. There are also support vector machines \cite{VapnikBook,SVM,bishop_pattern_2006,scholkopf_learning_2001}, which have been developed by the machine-learning community for classification tasks but are now used as surrogates in science and engineering as well, see, e.g., \cite{forrester_recent_2009,Basudhar2010,dribusch_multifidelity_2010,peherstorfer15MFMC}. 

\subsection{Outer-loop applications}\label{sec:intro-ol}
We focus on three outer-loop applications for which a range of multifidelity methods exist: uncertainty propagation, statistical inference, and optimization.

\paragraph{Uncertainty propagation} In uncertainty propagation, the model input is described by a random variable and one is interested in statistics of the model output. Using Monte Carlo simulation to estimate statistics of the model output often requires a large number of model evaluations to achieve accurate approximations of the statistics.
A multifidelity method that combines outputs from computationally cheap low-fidelity models with outputs from the high-fidelity model can lead to significant reductions in runtime and provide unbiased estimators of the statistics of the high-fidelity model outputs \cite{giles_multi-level_2008,ng_multifidelity_2014,narayan_stochastic_2014,teckentrup_multilevel_2014,Peherstorfer2016490}. Note that we consider probabilistic approaches to uncertainty propagation only; other approaches to uncertainty propagation are, e.g., fuzzy set approaches \cite{bernardini_what_1999} and worst-case scenario analysis \cite{Babuska2005DetUQ}.

\paragraph{Statistical inference}
In inverse problems, an indirect observation of a quantity of interest is given. A classical example is that limited and noisy observations of a system output are given and one wishes to estimate the input of the system.
In statistical inference \cite{tarantola_inverse_1982,TarantolaBook,KaipoBook,ANU:7701764}, the unknown input is modeled as a random variable and one is interested in sampling the distribution of this random variable to assess the uncertainty associated with the input estimation. Markov chain Monte Carlo (MCMC) methods provide one way to sample the distribution of the input random variable. MCMC is an outer-loop application that requires evaluating the high-fidelity model many times. Multifidelity methods in MCMC typically use multi-stage adaptive delayed acceptance formulations that leverage low-fidelity models to speed up the sampling \cite{christen_markov_2005,efendiev_preconditioning_2006,WRCR:WRCR13044,cui_data-driven_2014}.

\paragraph{Optimization} The goal of optimization is to find an input that leads to an optimal model output with respect to a given objective function. Optimization is typically solved using an iterative process that requires evaluations of the model in each iteration.
Multifidelity optimization reduces the runtime of the optimization process by using low-fidelity models to accelerate the search \cite{booker_rigorous_1999,keane_wing_2003,forrester_multi-fidelity_2007,forrester_recent_2009} or by using a low-fidelity model in conjunction with adaptive corrections and a trust-region model-management scheme \cite{alexandrov_trust-region_1998,alexandrov_approximation_2001,booker_rigorous_1999,march_provably_2012,Queipo20051}.
Other multifidelity optimization methods build a surrogate using evaluations from multiple models, and then optimize using this surrogate. For example, efficient global optimization (EGO) is a multifidelity optimization method that adaptively constructs a low-fidelity model by interpolating the objective function corresponding to the high-fidelity model with Gaussian process regression (kriging) \cite{jones_efficient_1998}.

\subsection{Outline of the paper}\label{sec:intro-outline}
The remainder of this paper focuses on model management strategies.
Section~\ref{sec:MM:strategies} overviews the model management strategies of adaptation, fusion, and filtering. Sections~\ref{sec:MM:UQ}--\ref{sec:MM:Optimization} survey specific techniques in the context of uncertainty propagation, inference, and optimization, respectively. The outlook in Section~\ref{sec:Outlook} closes the survey.

%% file: mm.tex
\section{Multifidelity model management strategies}
\label{sec:MM:strategies}
Model management in multifidelity methods defines how different models are employed during execution of the outer loop and how outputs from different models are combined. Models are managed such that low-fidelity models are leveraged for speedup, while judicious evaluations of the high-fidelity model establish accuracy and/or convergence of the outer-loop result. This section describes a categorization of model management methods into three types of strategies. The following sections then survey specific model management methods in the context of uncertainty propagation, statistical inference, and optimization.

As shown in Figure~\ref{fig:MMTypes}, we distinguish between three types of model management strategies: adaptation, fusion, and filtering.

\subsection{Adaptation} The first model management strategy uses adaptation to enhance the low-fidelity model with information from the high-fidelity model while the computation proceeds. One example of model management based on adaptation is global optimization with EGO, where a kriging model is adapted in each iteration of the optimization process \cite{jones_efficient_1998,venter_construction_1998}. Another example is the correction of low-fidelity model outputs via updates, which are derived from the high-fidelity model. It is common to use additive updates, which define the correction based on the difference between sampled high-fidelity and low-fidelity outputs, and/or multiplicative updates, which define the correction based on the ratio between sampled high-fidelity and low-fidelity outputs \cite{alexandrov_trust-region_1998,alexandrov_approximation_2001}. The correction model is then typically built using Taylor series expansion based on gradients, and possibly also on higher-order derivative information \cite{michael_eldred_second-order_2004}. In \cite{kennedy_predicting_2000}, low-fidelity models are corrected (calibrated) with Gaussian process models to best predict the output of the high-fidelity model.
Another multifidelity adaptation strategy is via adaptive model reduction, where projection-based reduced models are efficiently adapted as more data of the high-fidelity model become available during solution of the outer-loop application problem. Key to online adaptive model reduction is an efficient adaptation process. In \cite{Peherstorfer201521,aDEIM}, the basis and operators of projection-based reduced models are adapted with low-rank updates. In \cite{NME:NME4800}, an $h$-adaptive refinement of the basis vectors uses clustering algorithms to learn and adapt a reduced basis from high-fidelity model residuals. The work \cite{Amsallem2015} adapts localized reduced bases to smooth the transition from one localized reduced basis to another localized basis.

\subsection{Fusion} The second model management strategy is based on information fusion. Approaches based on fusion evaluate low- and high-fidelity models and then combine information from all outputs. An example from uncertainty propagation is the control variate framework \cite{Bratley,MCMethods,nelson_control_1987}, where the variance of Monte Carlo estimators is reduced by exploiting the correlation between high- and low-fidelity models. The control variate framework leverages a small number of high-fidelity model evaluations to obtain unbiased estimators of the statistics of interest, together with a large number of low-fidelity model evaluations to obtain an estimator with a low variance. Another example from uncertainty propagation is the fusion framework introduced in \cite{koutsourelakis_accurate_2009}, which is based on Bayesian regression.

Co-kriging is another example of a multifidelity method that uses model management based on fusion. Co-kriging derives a model from multiple information sources, e.g., a low- and a high-fidelity model \cite{annels_geostatistical_1991,myers_matrix_1982,Perdikaris20150018}. Co-kriging is often used in the context of optimization if gradient information of the high-fidelity model is available, see \cite{forrester_multi-fidelity_2007}. The work \cite{laurenceau_building_2008} compares kriging and co-kriging models on aerodynamic test functions. In \cite{yamazaki_derivative-enhanced_2012}, gradients are computed cheaply with the adjoint method and then used to derive a co-kriging model for design optimization in large design spaces. In \cite{Han2013177}, co-kriging with gradients and further developments of co-kriging are compared for approximating aerodynamic models of airfoils. 

\subsection{Filtering} The third model management strategy is based on filtering, where the high-fidelity model is invoked following the evaluation of a low-fidelity filter. This might entail evaluating the high-fidelity model only if the low-fidelity model is deemed inaccurate, or it might entail evaluating the high-fidelity model only if the candidate point meets some criterion based on the low-fidelity evaluation.
One example of a multifidelity filtering strategy is a multi-stage MCMC algorithm. For example, in two-stage MCMC \cite{christen_markov_2005,Fox97samplingconductivity}, a candidate sample needs to be first accepted by the likelihood induced by the low-fidelity model before the high-fidelity model is evaluated at the candidate sample. As another example, in the multifidelity stochastic collocation approach in \cite{narayan_stochastic_2014}, the stochastic space is explored with the low-fidelity model to derive sampling points at which the high-fidelity model is then evaluated. A third example is multifidelity importance sampling, where the sampling of the high-fidelity model is guided by an importance sampling biasing distribution that is constructed with a low-fidelity model \cite{Peherstorfer2016490}.

%% file: uq.tex
\section{Multifidelity model management in uncertainty propagation}
\label{sec:MM:UQ}
Inputs of models are often formulated as random variables to describe the stochasticity of the system of interest. With random inputs, the output of the model becomes a random variable as well. Uncertainty propagation aims to estimate statistics of the output random variable \cite{UQDef}. Sampling-based methods for uncertainty propagation evaluate the model at a large number of inputs and then estimate statistics from the corresponding model outputs. Examples of sampling-based methods are Monte Carlo simulation and stochastic collocation approaches. In this section, we review multifidelity approaches for sampling-based methods in uncertainty propagation. These multifidelity approaches shift many of the model evaluations to low-fidelity models while evaluating the high-fidelity model a small number of times to establish unbiased estimators. Section~\ref{sec:MM:UQ:MonteCarlo} introduces the problem setup and briefly overviews the Monte Carlo simulation method. Sections~\ref{sec:MM:UQ:ContVar}--\ref{sec:UQ:Other} discuss multifidelity methods for Monte Carlo based on control variates, importance sampling, and other techniques, respectively. Multifidelity methods for stochastic collocation are discussed in Section~\ref{sec:MM:UQ:StochColl}.

\subsection{Uncertainty propagation and Monte Carlo simulation}
\label{sec:MM:UQ:MonteCarlo}
Consider the high-fidelity model $\fhigh: \Zcal \to \Ycal$ and let the uncertainties in the inputs be represented by a random variable $Z$ with probability density function $p$. At this point, the only assumption we make on the random variable $Z$ is that the distribution is absolutely continuous such that a density function exists. In particular, the random variable $Z$ can be a non-Gaussian random variable. The goal of uncertainty propagation is to estimate statistics of the random variable $\fhigh(Z)$, e.g., the expectation
\begin{equation}
\mathbb{E}[\fhigh] = \int_{\Zcal} \fhigh(\dv)p(\dv)\mathrm d\dv\,,
\label{eq:MM:UQ:ExpVal}
\end{equation}
and the variance
\begin{equation}
\Var[\fhigh] = \mathbb{E}[\fhigh^2] - \mathbb{E}[\fhigh]^2\,,
\label{eq:MM:UQ:Var}
\end{equation}
which we assume to exist.

The Monte Carlo method draws $\nrpts \in \mathbb{N}$ independent and identically distributed (i.i.d.) realizations $\dv_1, \dots, \dv_{\nrpts} \in \Zcal$ of the random variable $Z$, and estimates the expectation $\mathbb{E}[\fhigh]$ as
\begin{equation}
\MCEst^{\hi}_{\nrpts} = \frac{1}{\nrpts}\sum_{i = 1}^{\nrpts} \fhigh(\dv_i)\,.
\label{eq:MM:UQ:MCEst}
\end{equation}
The Monte Carlo estimator is an unbiased estimator $\MCEst^{\hi}_{\nrpts}$ of $\mathbb{E}[\fhigh]$, which means that $\mathbb{E}[\MCEst^{\hi}_{\nrpts}] = \mathbb{E}[\fhigh]$. The mean squared error (MSE) of the Monte Carlo estimator $\MCEst^{\hi}_{\nrpts}$ therefore is
\begin{equation}
e(\MCEst^{\hi}_{\nrpts}) = \frac{\Var[\fhigh]}{\nrpts}\,.
\label{eq:MM:UQ:MCRMSE}
\end{equation}

The convergence rate $\mathcal{O}(\nrpts^{-1/2})$ of the root mean squared error (RMSE) $\sqrt{e(\MCEst^{\hi}_{\nrpts})}$ is low if compared to deterministic quadrature rules, see Section~\ref{sec:MM:UQ:StochColl}; however, the rate is independent of the smoothness of the integrand and the dimension $d$ of the input $\dv$, which means that the Monte Carlo method is well-suited for high dimensions $d$, and, in fact, is often the only choice available if $d$ is large. Typically more important in practice, however, is the pre-asymptotic behavior of the RMSE of the Monte Carlo estimator. In the pre-asymptotic regime, the variance $\Var[\fhigh]$ dominates the RMSE. Variance reduction techniques reformulate the estimation problem such that a function with a lower variance is integrated instead of directly integrating $\fhigh(Z)$. Examples of variance reduction techniques are antithetic variates, control variates, importance sampling, conditional Monte Carlo sampling, and stratified sampling \cite{MCMethods,RobertBook}. Variance reduction techniques often exploit the correlation between the random variable $\fhigh(Z)$ of interest and an auxiliary random variable. Multifidelity methods construct the auxiliary random variable using low-fidelity models. We discuss multifidelity methods for variance reduction based on control variates in Section~\ref{sec:MM:UQ:ContVar}, and variance reduction based on importance sampling in Section~\ref{sec:MM:UQ:ImportanceSampling}.

\subsection{Multifidelity uncertainty propagation based on control variates}
\label{sec:MM:UQ:ContVar}
The control variate framework \cite{MCMethods,Bratley,nelson_control_1987} aims to reduce the estimator variance of a random variable by exploiting the correlation with an auxiliary random variable. In the classical control variate method, as discussed in, e.g., \cite{MCMethods}, the statistics of the auxiliary random variable is known. Extensions relax this requirement by estimating the statistics of the auxiliary random variable from prior information \cite{emsermann_improving_2002,pasupathy_control-variate_2012}. We now discuss multifidelity approaches that construct auxiliary random variables from low-fidelity models.

\subsubsection{Control variates based on low-fidelity models}
Consider the high-fidelity model $\fhigh$ and $k \in \mathbb{N}$ low-fidelity models $\fl^{(1)}, \dots, \fl^{(k)}$. In \cite{ng_multifidelity_2014,peherstorfer15MFMC}, a multifidelity method is introduced that uses the random variables $\fl^{(1)}(Z), \dots, \fl^{(k)}(Z)$ stemming from the low-fidelity models as control variates for estimating statistics of the random variable $\fhigh(Z)$ of the high-fidelity model. An optimal model management is derived that minimizes the MSE of the multifidelity estimator for a given computational budget. In the numerical experiments, high-fidelity finite element models are combined with projection-based models, data-fit models, and support vector machines, which demonstrates that the multifidelity approach is applicable to a wide range of low-fidelity model types.

Let $\nrpts_0 \in \mathbb{N}$ be the number of high-fidelity model evaluations. Let $\nrpts_i \in \mathbb{N}$ be the number of evaluations of the low-fidelity model $\fl^{(i)}$ for $i = 1, \dots, k$, where $0 < \nrpts_0 \leq \nrpts_1 \leq \dots \leq \nrpts_k$. The multifidelity approach presented in \cite{peherstorfer15MFMC} draws $\nrpts_k$ realizations
\begin{equation}
\dv_1, \dots, \dv_{\nrpts_k}
\label{eq:MM:UQ:Realizations}
\end{equation}
from the random variable $Z$ and computes the model outputs $\fhigh(\dv_1), \dots, \fhigh(\dv_{\nrpts_0})$ and
\begin{equation}
\fl^{(i)}(\dv_1), \dots, \fl^{(i)}(\dv_{\nrpts_i})
\label{eq:MM:UQ:ModelOutputs}
\end{equation}
for $i = 1, \dots, k$. These model outputs are used to derive Monte Carlo estimates
\begin{equation}
\MCEst^{\hi}_{\nrpts_0} = \frac{1}{\nrpts_0}\sum_{j = 1}^{\nrpts_0}\fhigh(\dv_j)\,,\qquad\MCEst^{(i)}_{\nrpts_i} = \frac{1}{\nrpts_i}\sum_{j = 1}^{\nrpts_i}\fl^{(i)}(\dv_j)\,,\qquad i = 1, \dots, k\,,
\label{eq:MM:UQ:MCEstOne}
\end{equation}
and
\begin{equation}
\MCEst^{(i)}_{\nrpts_{i - 1}} = \frac{1}{\nrpts_{i - 1}}\sum_{j = 1}^{\nrpts_{i - 1}}\fl^{(i)}(\dv_j)\,,\qquad i = 1, \dots, k\,.
\label{eq:MM:UQ:MCEstTwo}
\end{equation}
Note that the estimates \eqref{eq:MM:UQ:MCEstTwo} use the first $\fl^{(i)}(\dv_1), \dots, \fl^{(i)}(\dv_{\nrpts_{i - 1}})$ model outputs of \eqref{eq:MM:UQ:ModelOutputs} only, whereas the estimate \eqref{eq:MM:UQ:MCEstOne} uses all $\nrpts_i$ model outputs \eqref{eq:MM:UQ:ModelOutputs} for $i = 1, \dots, k$. Following \cite{peherstorfer15MFMC}, the multifidelity estimator of $\mathbb{E}[\fhigh]$ is
\begin{equation}
\MCEst^{\text{MF}} = \MCEst^{\hi}_{\nrpts_0} + \sum_{i = 1}^k\alpha_i\left(\MCEst^{(i)}_{\nrpts_i} - \MCEst^{(i)}_{\nrpts_{i - 1}}\right)\,.
\label{eq:MM:UQ:ContVarEst}
\end{equation}
The control variate coefficients $\alpha_1, \dots, \alpha_k \in \mathbb{R}$ balance the term $\MCEst^{\hi}_{\nrpts_0}$ stemming from the high-fidelity model and the terms $\MCEst^{(i)}_{\nrpts_i} - \MCEst^{(i)}_{\nrpts_{i - 1}}$ from the low-fidelity models. The multifidelity estimator \eqref{eq:MM:UQ:ContVarEst} based on the control variate framework evaluates the high- and the low-fidelity model and fuses both outputs into an estimate of the statistics of the high-fidelity model. The multifidelity estimator \eqref{eq:MM:UQ:ContVarEst} therefore uses a model management based on fusion, see Section~\ref{sec:MM:strategies}. We note that \eqref{eq:MM:UQ:ContVarEst} could also be viewed as a correction, although the correction is to the estimators stemming from the low-fidelity models, not to the low-fidelity model outputs directly.

\paragraph{Properties of the multifidelity estimator}
The multifidelity estimator $\MCEst^{\text{MF}}$ is an unbiased estimator of $\mathbb{E}[\fhigh]$ because
\[
\mathbb{E}[\MCEst^{\text{MF}}] = \mathbb{E}[\MCEst^{\hi}_{\nrpts_0}] + \sum_{i = 1}^k\alpha_i\mathbb{E}[\MCEst^{(i)}_{\nrpts_i} - \MCEst^{(i)}_{\nrpts_{i - 1}}] = \mathbb{E}[\fhigh]\,.
\]
Therefore, the MSE of the estimator $\MCEst^{\text{MF}}$ is equal to the variance $\Var[\MCEst^{\text{MF}}]$ of the estimator, $e(\MCEst^{\text{MF}}) = \Var[\MCEst^{\text{MF}}]$. The costs of the multifidelity estimator are
\[
c(\MCEst^{\text{MF}}) = \nrpts_0\chigh + \sum_{i = 1}^k\nrpts_i\clow^{(i)} = \boldsymbol{\nrpts}^T\bfc\,,
\]
where $\boldsymbol{\nrpts} = [\nrpts_0, \nrpts_1, \dots, \nrpts_k]^T$ and $\bfc = [\chigh, \clow^{(1)}, \dots, \clow^{(k)}]^T$. The high-fidelity model is evaluated at $\nrpts_0$ realizations and the low-fidelity model $\fl^{(i)}$ at $\nrpts_i$ realizations of $Z$, for $i = 1, \dots, k$.

\subsubsection{Multifidelity Monte Carlo}
The multifidelity estimator \eqref{eq:MM:UQ:ContVarEst} depends on the control variate coefficients $\alpha_1, \dots, \alpha_k$ and on the number of model evaluations $\nrpts_0, \nrpts_1, \dots, \nrpts_k$. In \cite{ng_multifidelity_2014,peherstorfer15MFMC}, these parameters are chosen such that the MSE of the estimator \eqref{eq:MM:UQ:ContVarEst} is minimized for a given computational budget $\gamma \in \mathbb{R}_{+}$. The solution to the optimization problem
\begin{equation}
\begin{aligned}
\min_{\substack{\nrpts_0, \nrpts_1, \dots, \nrpts_k\\\alpha_1, \dots, \alpha_k}} \quad & e(\MCEst^{\text{MF}})\\
\text{s.t.} \quad & \nrpts_0 > 0\\
& \nrpts_i \geq \nrpts_{i - 1}\,,\qquad i = 1, \dots, k\\
& \boldsymbol{\nrpts}^T\bfc = \gamma
\end{aligned}
\label{eq:MM:UQ:ContVarOptiProblem}
\end{equation}
gives the coefficients $\alpha^*_1, \dots, \alpha_k^*$ and the number of model evaluations $\nrpts^*_0, \dots, \nrpts^*_k$ that minimize the MSE of the multifidelity estimator $\MCEst^{\text{MF}}$ for the given computational budget $\gamma$. The constraints impose that $0 < \nrpts_0^* \leq \nrpts_1^* \leq \dots \leq \nrpts_k^*$ and that the costs $c(\MCEst^{\text{MF}})$ of the estimator equal the computational budget $\gamma$.

\paragraph{Variance of the multifidelity estimator}
Since the multifidelity estimator $\MCEst^{\text{MF}}$ is unbiased, we have $e(\MCEst^{\text{MF}}) = \Var[\MCEst^{\text{MF}}]$, and therefore the objective of minimizing the MSE $e(\MCEst^{\text{MF}})$ can be replaced with the variance $\Var[\MCEst^{\text{MF}}]$ in the optimization problem \eqref{eq:MM:UQ:ContVarOptiProblem}. The variance $\Var[\MCEst^{\text{MF}}]$ of the multifidelity estimator $\MCEst^{\text{MF}}$ is
\begin{equation}
\Var[\MCEst^{\text{MF}}] = \frac{\sigma_{\hi}^2}{\nrpts_0} + \sum_{i = 1}^k\left(\frac{1}{\nrpts_{i - 1}} - \frac{1}{\nrpts_i}\right)\left(\alpha^2_i\sigma_i^2 - 2\alpha_i\rho_i\sigma_{\hi}\sigma_i\right)\,,
\label{eq:MM:UQ:ContVarVariance}
\end{equation}
where $-1 \leq \rho_i \leq 1$ is the Pearson correlation coefficient of the random variable $\fhigh(Z)$ and $\fl^{(i)}(Z)$ for $i = 1, \dots, k$. The quantities
\[
\sigma_{\hi}^2 = \Var[\fhigh]\,,\qquad \sigma_i^2 = \Var[\fl^{(i)}]\,,\qquad i = 1, \dots, k\,,
\]
are the variances of $\fhigh(Z)$ and $\fl^{(i)}(Z)$, respectively.

\paragraph{Optimal selection of the number of samples and control variate coefficients}
Under certain conditions on the low- and the high-fidelity model, the optimization problem \eqref{eq:MM:UQ:ContVarOptiProblem} has a unique, analytic solution \cite{peherstorfer15MFMC}. The optimal control variate coefficients are
\begin{equation}
\alpha^*_i = \rho_i\frac{\sigma_{\hi}}{\sigma_i}\,,\qquad i = 1, \dots, k\,.
\label{eq:MM:UQ:ContVarOptiParameters}
\end{equation}
The optimal number of evaluations $\nrpts_0^*, \nrpts_1^*, \dots, \nrpts_k^*$ are
\begin{equation}
  \nrpts_0^* = \frac{\gamma}{\bfc^T\bfr}\,,\qquad \nrpts_i^* = r_i\nrpts_0\,,\qquad i = 1, \dots, k\,,
  \label{eq:MM:UQ:OptimalNrEvals}
\end{equation}
where the components of the vector $\bfr = [1, r_1, \dots, r_k]^T \in \mathbb{R}^{k + 1}$ are given as
\begin{equation}
  r_i = \sqrt{\frac{\chigh(\rho_i^2 - \rho_{i + 1}^2)}{\clow^{(i)}(1 - \rho_1^2)}}\,,\qquad i = 1, \dots, k\,.
  \label{eq:MM:UQ:OptimalR}
\end{equation}
Note that the convention $\rho_{k + 1} = 0$ is used in \eqref{eq:MM:UQ:OptimalR}. We refer to \cite{peherstorfer15MFMC} for details.

\paragraph{Interaction of models in multifidelity Monte Carlo}
We compare the multifidelity estimator $\MCEst^{\text{MF}}$ to a benchmark Monte Carlo estimator $\MCEst^{\text{MC}}$ that uses the high-fidelity model alone. The multifidelity estimator $\MCEst^{\text{MF}}$ and the benchmark estimator $\MCEst^{\text{MC}}$ have the same costs $\gamma$. With the MSE $e(\MCEst^{\text{MF}})$ of the multifidelity estimator and the MSE $e(\MCEst^{\text{MC}})$ of the benchmark Monte Carlo estimator, the variance reduction ratio is
\begin{equation}
\frac{e(\MCEst^{\text{MF}})}{e(\MCEst^{\text{MC}})} = \left(\sqrt{1 - \rho_1^2} + \sum_{i = 1}^k\sqrt{\frac{\clow^{(i)}}{\chigh}(\rho_i^2 - \rho_{i + 1}^2)}\right)^2\,.
\label{eq:MM:UQ:VarRedRatio}
\end{equation}
The ratio \eqref{eq:MM:UQ:VarRedRatio} quantifies the variance reduction achieved by the multifidelity estimator compared to the benchmark Monte Carlo estimator. The variance reduction ratio is a sum over the costs $\chigh, \clow^{(1)}, \dots, \clow^{(k)}$ and the correlation coefficients $\rho_1, \dots, \rho_k$ of all models in the multifidelity estimator. This shows that the contribution of a low-fidelity model to the variance reduction of the multifidelity estimator cannot be determined by the properties of that low-fidelity model alone but only by taking into account all other models that are used in the multifidelity estimator. Thus, the interaction between the models is what drives the efficiency of the multifidelity estimator $\MCEst^{\text{MF}}$. We refer to \cite{peherstorfer15MFMC} for an in-depth discussion of the interaction between the models and a more detailed analysis.

\paragraph{Efficiency of the multifidelity estimator}
It is shown in \cite{peherstorfer15MFMC} that the MFMC estimator $\MCEst^{\text{MF}}$ is computationally cheaper than the benchmark Monte Carlo estimator that uses the high-fidelity model $\fhigh$ alone if
\begin{equation}
\sqrt{1 - \rho_1^2} + \sum_{i = 1}^k\sqrt{\frac{\clow^{(i)}}{\chigh}(\rho_i^2 - \rho_{i + 1}^2)} < 1
\label{eq:MM:UQ:BetterThanMC}
\end{equation}
The inequality \eqref{eq:MM:UQ:BetterThanMC} emphasizes that both correlation and costs of the models are critical for an efficient multifidelity estimator.

\paragraph{Algorithm}
Algorithm~\ref{alg:MM:UQ:MFMC} summarizes the multifidelity Monte Carlo method as presented in \cite{peherstorfer15MFMC}. Inputs are the models $\fhigh, \fl^{(1)}, \dots, \fl^{(k)}$ and the variances $\sigma_{\text{hi}}$ and $\sigma_1, \dots, \sigma_k$. The inputs $\rho_i$ are the correlation coefficients of the random variable $\fhigh(Z)$ stemming from the high-fidelity model and the random variables $\fl^{(i)}(Z)$ for $i = 1, \dots, k$. The costs of the models are $\chigh, \clow^{(1)}, \dots, \clow^{(k)}$ and the computational budget is $\gamma$. Line~\ref{alg:MM:UQ:MFMC:Ordering} of Algorithm~\ref{alg:MM:UQ:MFMC} ensures that the ordering of the models is used that minimizes the MSE of the multifidelity estimator, see \cite[Section~3.5]{peherstorfer15MFMC} for details. Line~\ref{alg:MM:UQ:MFMC:R} defines the vector of ratios $\bfr = [r_0, r_1, \dots, r_k]^T$, cf.~\eqref{eq:MM:UQ:OptimalR}. The numbers of model evaluations $\nrpts_0^*, \nrpts_1^*, \dots, \nrpts_k^*$ of the models $\fhigh, \fl^{(1)}, \dots, \fl^{(k)}$ are derived from $\bfr$ as in \eqref{eq:MM:UQ:OptimalNrEvals}. The control variate coefficients $\alpha_1^*, \dots, \alpha_k^*$ are obtained as in \eqref{eq:MM:UQ:ContVarOptiParameters}. In line~\ref{alg:MM:UQ:MFMC:Realizations}, $\nrpts_k^*$ realizations $\dv_1, \dots, \dv_{\nrpts_k^*}$ are drawn from the random variable $Z$. The high-fidelity model $\fhigh$ is evaluated at the realizations $\dv_1, \dots, \dv_{\nrpts_0^*}$ and models $\fl^{(i)}$ are evaluated at $\dv_1, \dots, \dv_{\nrpts_i^*}$ for $i = 1, \dots, k$. The multifidelity estimate is obtained as in \eqref{eq:MM:UQ:ContVarEst} and returned.

\begin{algorithm}[t]
\caption{Multifidelity Monte Carlo}\label{alg:MM:UQ:MFMC}
\begin{algorithmic}[1]
\Procedure{MFMC}{$\fhigh, \fl^{(1)}, \dots, \fl^{(k)}, \sigma_{\text{hi}}, \sigma_1, \dots, \sigma_k, \rho_{1}, \dots, \rho_{k}, \chigh, \clow^{(1)}, \dots, \clow^{(k)}, \gamma$}
\State Ensure $\fhigh, \fl^{(1)}, \dots, \fl^{(k)}$ are ordered as described in \cite[Section~3.5]{peherstorfer15MFMC}\label{alg:MM:UQ:MFMC:Ordering}
\State Set $\rho_{k + 1} = 0$ and define vector $\bfr = [1, r_1, \dots, r_k]^T \in \mathbb{R}^{k + 1}_+$ as
\[
r_i = \sqrt{\frac{\chigh(\rho_{i}^2 - \rho^2_{i + 1})}{\clow^{(i)}(1 - \rho_{1}^2)}}\,,\qquad i = 1, \dots, k
\]\label{alg:MM:UQ:MFMC:R}
\State Select number of model evaluations $\boldsymbol{\nrpts^*} \in \mathbb{R}^{k + 1}_+$ as
\[
\boldsymbol{\nrpts^*} = \left[\frac{\gamma}{\boldsymbol{c}^T\bfr}, r_1\nrpts_0^*, \dots, r_k\nrpts_0^*\right]^T \in \mathbb{R}^{k + 1}_+
\]
\State Set coefficients $\boldsymbol{\alpha^*} = [\alpha_1^*, \dots, \alpha_k^*]^T \in \mathbb{R}^k$ to
\[
\alpha_i^* = \frac{\rho_{i}\sigma_{\text{hi}}}{\sigma_i}\,,\qquad i = 1, \dots, k
\]
\State Draw $\dv_1, \dots, \dv_{\nrpts_k^*} \in \Zcal$ realizations of $Z$\label{alg:MM:UQ:MFMC:Realizations}
\State Evaluate high-fidelity model $\fhigh$ at realizations $\dv_1, \dots, \dv_{\nrpts_0^*}$
\State Evaluate model $\fl^{(i)}$ at realizations $\dv_1, \dots, \dv_{\nrpts_i^*}$ for $i = 1, \dots, k$\label{alg:MM:UQ:MFMC:ModelEvals}
\State Compute the multifidelity estimate $\MCEst^{\text{MF}}$ as in \eqref{eq:MM:UQ:ContVarEst}
\State \Return multifidelity estimate $\MCEst^{\text{MF}}$
\EndProcedure
\end{algorithmic}
\end{algorithm}

\subsubsection{Other uses of control variates as a multifidelity technique} In \cite{boyaval_fast_2012,boyaval_variance_2010}, the reduced basis method is used to construct control variates. The reduced basis models are built with greedy algorithms that use \emph{a posteriori} error estimators to particularly target variance reduction. The work \cite{VidalCodina2015700} uses error estimators to combine reduced basis models with control variates. The StackMC method presented in \cite{tracey_using_2013} successively constructs machine-learning-based low-fidelity models and combines them with the control variate framework. In \cite{ng_monte_2015}, the multifidelity control variate method is used in the context of optimization, where information of previous iterations of the optimization problem are used as control variate. This means that data from previous iterations serve as a kind of low-fidelity ``model''.

The multilevel Monte Carlo method \cite{heinrich_multilevel_2001,giles_multi-level_2008} uses the control variate framework to combine multiple low-fidelity models with a high-fidelity model. Typically, in multilevel Monte Carlo, the low-fidelity models are coarse-grid approximations, where the accuracy and costs can be controlled by a discretization parameter. The properties of the low-fidelity models are therefore often described with rates. For example, the rate of the decay of the variance of the difference of two successive coarse-grid approximations and the rate of the increase of the costs with finer grids play a critical role in determining the efficiency of the multilevel Monte Carlo method. Additionally, rates are used to determine the number of evaluations of each low-fidelity model and the high-fidelity model, see, e.g., \cite[Theorem 1]{cliffe_multilevel_2011}.
In the setting of stochastic differential equations and coarse-grid approximations, multilevel Monte Carlo has been very successful, see, e.g., \cite{cliffe_multilevel_2011,teckentrup_further_2013}, the recent advances on multi-index Monte Carlo \cite{haji-ali_multi-index_2015}, and the nesting of multilevel Monte Carlo and control variates \cite{nobile_multi_2015} for detailed studies and further references. 

\subsection{Multifidelity uncertainty propagation based on importance sampling}
\label{sec:MM:UQ:ImportanceSampling}
Importance sampling \cite{MCMethods} uses a problem-dependent sampling strategy. The goal is an estimator with a lower variance than a Monte Carlo estimator such as \eqref{eq:MM:UQ:MCEst}. The problem-dependent sampling means that samples are drawn from a biasing distribution, instead of directly from the distribution of the random variable $Z$ of interest, and then the change of the distribution is compensated with a re-weighting. Importance sampling is particularly useful in the case of rare event simulation, where the probability of the event of interest is small and therefore many realizations of the random variable $Z$ are necessary to obtain a Monte Carlo estimate of reasonable accuracy. Importance sampling with a suitable biasing distribution can explicitly target the rare event and reduce the number of realizations required to achieve an acceptable accuracy. The challenge of importance sampling is the construction of a biasing distribution, which usually is problem-dependent and typically requires model evaluations. We discuss multifidelity methods that use low-fidelity models for the construction of biasing distributions.

\subsubsection{Importance sampling}
Consider the indicator function $I_{\hi}: \Zcal \to \{0, 1\}$ defined as
\[
I_{\hi}(\bfz) = \begin{cases}
1\,,\qquad & \fhigh(\bfz) < 0\,,\\
0\,,\qquad & \fhigh(\bfz) \geq 0
\end{cases}\,.
\]
We define the set $\mathcal{I} = \{\dv \in \Zcal \,|\, I_{\hi}(\bfz) = 1\}$. The goal is to estimate the probability of the event $Z^{-1}(\mathcal{I})$, which is $\mathbb{E}_p[I_{\hi}]$, with importance sampling. Note that we now explicitly denote in the subscript of $\mathbb{E}$ with respect to which distribution the expectation is taken.

\paragraph{Step 1: Construction of biasing distribution}
Traditionally, importance sampling consists of two steps. In the first step, the biasing distribution with density $q$ is constructed. Let $Z^{\prime}$ be the biasing random variable with the biasing density $q$. Recall that the input random variable $Z$ with the nominal distribution has the nominal density $p$. Let
\[
\operatorname{supp}(p) = \{\dv \in \Zcal: p(\dv) > 0\}
\]
be the support of the density $p$. If the support of the nominal density $p$ is a subset of the support of the biasing density $q$, i.e., $\operatorname{supp}(p) \subset \operatorname{supp}(q)$, then the expectation $\mathbb{E}_p[I_{\hi}]$ can be rewritten in terms of the biasing density $q$ as
\[
\mathbb{E}_p[I_{\hi}] = \int_{\Zcal} I_{\hi}(\dv)p(\dv)\mathrm d\dv = \int_{\Zcal} I_{\hi}(\dv^{\prime})q(\dv^{\prime})\frac{p(\dv^{\prime})}{q(\dv^{\prime})}\mathrm d\dv^{\prime} = \mathbb{E}_q\left[I_{\hi}\frac{p}{q}\right]\,,
\]
where the ratio $p/q$ serves as a weight.

\paragraph{Step 2: Deriving an importance sampling estimate}
In the second step, the importance sampling estimator
\begin{equation}
\MCEst^{\text{IS}}_{\nrpts} = \frac{1}{\nrpts}\sum_{i = 1}^{\nrpts} I_{\hi}(\dv^{\prime}_i)\frac{p(\dv^{\prime}_i)}{q(\dv^{\prime}_i)}
\label{eq:MM:UQ:ISEst}
\end{equation}
is evaluated for realizations $\dv^{\prime}_1, \dots, \dv^{\prime}_{\nrpts} \in \Zcal$ of the random variable $Z^{\prime}$. The MSE of the estimator \eqref{eq:MM:UQ:ISEst} is
\begin{equation}
e(\MCEst^{\text{IS}}_{\nrpts}) = \frac{\Var_q[I_{\hi}\frac{p}{q}]}{\nrpts}\,.
\label{eq:MM:UQ:ISRMSE}
\end{equation}

\paragraph{Variance of importance sampling estimator}
The variance in \eqref{eq:MM:UQ:ISRMSE} is with respect to the biasing density $q$, cf.~Section~\ref{sec:MM:UQ:MonteCarlo} and the MSE of the Monte Carlo estimator \eqref{eq:MM:UQ:MCRMSE}. Therefore, the goal is to construct a biasing distribution with
\[
\Var_q\left[I_{\hi}\frac{p}{q}\right] < \Var_p[I_{\hi}]\,,
\]
to obtain an importance sampling estimator $\MCEst^{\text{IS}}_{\nrpts}$ that has a lower MSE than the Monte Carlo estimator $\MCEst^{\hi}_{\nrpts}$ for the same number of realizations $\nrpts$.

\subsubsection{Construction of the biasing distribution with low-fidelity models}
The multifidelity importance sampling approach introduced in \cite{Peherstorfer2016490} uses a low-fidelity model to construct the biasing distribution in the first step of importance sampling, and derives the statistics using high-fidelity model evaluations in step two. In that sense, multifidelity importance sampling uses a model management strategy based on filtering, see Section~\ref{sec:MM:strategies}.

In step one, the low-fidelity model $\fl$ is evaluated at a large number $n \in \mathbb{N}$ of realizations $\dv_1, \dots, \dv_{n}$ of the input random variable $Z$. This is computationally feasible because the low-fidelity model is cheap to evaluate. A mixture model $q$ of Gaussian distributions is fitted with the expectation-maximization algorithm to the set of realizations
\[
\{\dv_i \,|\, I_{\lo}(\dv_i) = 1\,, i = 1, \dots, n\}
\]
for which the low-fidelity model predicts the event of interest with the indicator function $I_{\lo}: \Zcal \to \{0, 1\}$
\[
I_{\lo}(\dv) = \begin{cases}
1\,,\qquad & \fl(\dv) < 0\,,\\
0\,,\qquad & \fl(\dv) \geq 0
\end{cases}\,.
\]
The mixture model $q$ serves as a biasing distribution. Note that other density estimation methods can be used instead of fitting a mixture model of Gaussian distributions with the expectation-maximization algorithm \cite{SilvermanDensity,libagf,SGDE}.

In step two, the high-fidelity model is evaluated at realizations $\dv_1^{\prime}, \dots, \dv_{\nrpts}^{\prime}$ from the biasing random variable $Z^{\prime}$ with biasing density $q$. From the high-fidelity model evaluations $\fhigh(\dv_1^{\prime}), \dots, \fhigh(\dv_{\nrpts}^{\prime})$ an estimate of the event probability $\mathbb{E}_p[I_{\hi}]$ is obtained. Under the condition that the support of the biasing density $q$ includes the support of the nominal density $p$, the multifidelity importance sampling approach leads to an unbiased estimator of the probability of the event. If the low-fidelity model is sufficiently accurate, then significant runtime savings can be obtained during the construction of the biasing distribution.  Note that using $I_{\lo}$ in the second step of the multifidelity approach would lead to a biased estimator of $\mathbb{E}_p[I_{\hi}]$ because $\fl$ is only an approximation of $\fhigh$ and thus $\mathbb{E}_p[I_{\hi}] \neq \mathbb{E}_p[I_{\lo}]$ in general.

\subsection{Other model management strategies for probability estimation and limit state function evaluation}
\label{sec:UQ:Other}
The multifidelity approaches for estimating $\mathbb{E}_p[I_{\hi}]$ that are discussed in Section~\ref{sec:MM:UQ:ImportanceSampling} first use the low-fidelity model to construct a biasing density $q$ and then the high-fidelity model to estimate the failure probability $\mathbb{E}_p[I_{\hi}]$ with importance sampling. Under mild conditions on the biasing density $q$ derived from the low-fidelity model, the importance sampling estimator using the high-fidelity model is an unbiased estimator of the failure probability $\mathbb{E}_p[I_{\hi}]$. In this section, we review multifidelity methods that combine low- and high-fidelity model evaluations to obtain an indicator function $\tilde{I}$ that approximates $I_{\hi}$ and that is computationally cheaper to evaluate than $I_{\hi}$. Thus, instead of exploiting the two-step importance sampling procedure as in Section~\ref{sec:MM:UQ:ImportanceSampling}, the techniques in this section leverage the low- and high-fidelity model to approximate $I_{\hi}$ with $\tilde{I}$ with the aim that the error $|\mathbb{E}_p[\tilde{I}] - \mathbb{E}_p[I_{\hi}]|$ is small.

The multifidelity approach introduced in \cite{li_evaluation_2010} is based on filtering to combine low- and high-fidelity model outputs to obtain an approximation $\tilde{I}$ of $I_{\hi}$. Let $\fl$ be a low-fidelity model and $\fhigh$ the high-fidelity model. Let further $\gamma > 0$ be a positive threshold value. The approach in \cite{li_evaluation_2010} considers the indicator function
\[
\tilde{I}(\bfz) = \begin{cases}
1\,, & \fl(\bfz) < -\gamma \text{ or } \big(|\fl(\bfz)| \leq \gamma \text{ and } \fhigh(\bfz) < 0\big)\\
0\,, & \text{else}
\end{cases}\,,
\]
which evaluates to 1 at an input $\bfz$ if either $\fl(\bfz) < -\gamma$ or $|\fl(\bfz)| \leq \gamma$ and $\fhigh(\bfz) < 0$. Evaluating the indicator function $\tilde{I}$ at $\bfz$ means that first the low-fidelity model $\fl$ is evaluated at $\bfz$. If $\fl(\bfz) < -\gamma$, then 1 is returned, and no high-fidelity model evaluation is necessary. Similarly, if $\fl(\bfz) > \gamma$, then the indicator function $\tilde{I}$ evaluates to 0 without requiring a high-fidelity model evaluation. However, if $|\fl(\bfz)| \leq \gamma$, then the input $\bfz$ lies near the failure boundary, and the high-fidelity model is evaluated to decide whether the indicator function returns 0 or 1. How often the high-fidelity model is evaluated is determined by the positive threshold value $\gamma$. The choice of $\gamma$ directly depends on the error of the low-fidelity model $\fl$ in a certain norm. If the error of $\fl$ is known, the work \cite{li_evaluation_2010} establishes a convergence theory under mild conditions on the error of $\fl$. In particular, the authors of \cite{li_evaluation_2010} show that the error $|\mathbb{E}_p[\tilde{I}] - \mathbb{E}_p[I_{\hi}]|$ can be reduced below any $\epsilon > 0$ by a choice of $\gamma$ that depends on the error of $\fl$. If the error of the low-fidelity model is unknown, the work \cite{li_evaluation_2010} proposes an iterative heuristic approach that avoids the choice of $\gamma$. In \cite{li_efficient_2011}, the multifidelity approach of \cite{li_evaluation_2010} is extended to importance sampling with the cross-entropy method. Similarly to the approaches in \cite{li_evaluation_2010,li_efficient_2011}, the work \cite{chen_accurate_2013} switches between low- and high-fidelity model evaluations by relying on \emph{a posteriori} error estimators for reduced basis models to decide if either the reduced model or the high-fidelity model should be used.

The multifidelity methods in \cite{Basudhar2010,dribusch_multifidelity_2010,aDEIM} all use a model management strategy based on adaptation to estimate a failure boundary or failure probability. In \cite{Basudhar2010}, a support vector machine is used to derive a low-fidelity model of the limit state function in failure probability estimation and design. The authors decompose the input space using decision boundaries obtained via the support vector machine, and so handle the discontinuities that arise when approximating the limit state function. An adaptive sampling scheme is introduced that refines the low-fidelity model along the failure boundary. 
In \cite{dribusch_multifidelity_2010}, another approach is introduced that uses a support vector machine to approximate the failure boundary. It is proposed to train the support vector machine with data obtained from a low- and a high-fidelity model. With the low-fidelity model, an initial approximation of the failure boundary is obtained by extensively sampling the input domain, which is computationally tractable because the low-fidelity model is cheap to evaluate. The training of the support vector machine for approximating the failure boundary corresponding to the high-fidelity model is then initialized with the approximate boundary obtained with the low-fidelity model. Additional samples are drawn in regions where the low- and the high-fidelity failure boundary differ, to refine the approximation of the high-fidelity failure boundary. The work \cite{aDEIM} presents an online adaptive reduced modeling approach, which is demonstrated in the context of failure probability estimation. To increase the accuracy of the reduced model, it is adapted to the failure boundary as the estimation of the failure probability proceeds. The adaptation is performed via low-rank updates to the basis matrix of the reduced model. The low-rank updates are derived from sparse samples of the high-fidelity model.

\subsection{Stochastic collocation and multifidelity}
\label{sec:MM:UQ:StochColl}
Stochastic collocation methods \cite{doi:10.1137/050645142,doi:10.1137/060663660,ANU:9260768} compute statistics such as the expectation \eqref{eq:MM:UQ:ExpVal} and the variance \eqref{eq:MM:UQ:Var} by using a deterministic quadrature rule instead of the Monte Carlo method. The quadrature rules are often based on sparse grids \cite{bungartz_sparse_2004,doi:10.1137/060663660} to perform the quadrature efficiently for high-dimensional inputs.

In \cite{Eldred2016}, statistics are computed using stochastic collocation, where the outputs of a low-fidelity model are corrected with a discrepancy model that accounts for the difference between the high- and the low-fidelity model. The discrepancy model is then used to derive either an additive correction, a multiplicative correction, or a weighted combination of additive and multiplicative corrections to the low-fidelity model outputs. Thus, this is another example of model management based on adaptation, see Section~\ref{sec:MM:strategies}. The authors of \cite{Eldred2016} point out that an adaptive refinement of the discrepancy model is necessary because the complexity of the discrepancy between the high- and the low-fidelity model varies distinctly in the stochastic domain. This is because low-fidelity models tend to approximate the high-fidelity model well only in certain regions of the stochastic domain, whereas in other regions they hardly match the high-fidelity model at all. 

Another multifidelity stochastic collocation method is presented in \cite{narayan_stochastic_2014}. This method is based on a filtering model management strategy. The low-fidelity model is first evaluated at a large number of collocation points to sample the stochastic domain. From these samples, a small number of points is selected via a greedy procedure, and the high-fidelity model is evaluated at these points. The state solutions of the high-fidelity model at the selected collocation points span a space in which approximations of the high-fidelity model states at all other sampling points are derived.

In \cite{teckentrup_multilevel_2014}, a multilevel stochastic collocation method uses a hierarchy of models to accelerate convergence. Low-fidelity models are coarse-grid approximations of the high-fidelity model. Similarly to the multilevel Monte Carlo method, a reduction of the computational complexity can be shown if the errors of the models in the hierarchy decay with a higher rate than the rate of the increase of the costs. We categorize this multilevel stochastic collocation method as model management based on fusion, because low- and high-fidelity model outputs are fused to derive an estimate of the statistics of the high-fidelity model.

%% file: invprob.tex
\section{Multifidelity model management in statistical inference}
\label{sec:MM:InvProb}
In a Bayesian setting, inverse problems are cast in a statistical formulation where the unknown input is modeled as a random variable and is described by its posterior distribution \cite{tarantola_inverse_1982,TarantolaBook,KaipoBook,ANU:7701764}.  MCMC is a popular way to sample from the posterior distribution. Statistical inference raises several computational challenges, including the design of MCMC sampling schemes, the construction of approximate models that can reduce the costs of MCMC sampling, and the development of alternatives to MCMC sampling such as variational approaches \cite{UQDef}. Detailed discussions of these and many other important aspects of statistical inference can be found in the literature, e.g., \cite{RobertBook,KaipoBook,ANU:7701764,MCBookLiu}. We focus here on a few specific aspects of statistical inference in which multifidelity methods have been used. In particular, we survey multifidelity methods that use a two-stage formulation of MCMC, where a candidate sample has to be first accepted by a low-fidelity model before it is passed on to be either accepted or rejected by the high-fidelity model. Section~\ref{sec:MM:InvProb:Bayesian} describes the problem setup. Section~\ref{sec:MM:InvProb:TwoStageMCMC} describes a two-stage MCMC framework and Section~\ref{sec:MM:InvProb:Adaptive} discusses a framework where a low-fidelity model is adapted while it is evaluated in an MCMC algorithm. A Bayesian approach to model and correct the error of low-fidelity models is discussed in Section~\ref{sec:MM:InvProb:Error}.

\subsection{Bayesian framework for inference}
\label{sec:MM:InvProb:Bayesian}
Consider the high-fidelity model $\fhigh$ that maps inputs $\dv \in \Zcal$ onto outputs $\ov \in \Ycal$. Let $p_0$ be a prior density that describes the input $\dv$ before any measurements. Let further $\obs$ be noisy observational data with the stochastic relationship
\begin{equation}
\obs = \fhigh(\dv) + \bfeps\,,
\label{eq:MM:InvProb:Bayesian:InvProb}
\end{equation}
where $\bfeps$ is a random vector that captures the measurement error, noise, and other uncertainties of the observation $\obs$. In the following, the random vector $\bfeps$ is modeled as a zero-mean Gaussian with covariance $\bfSigma_{\bfeps} \in \mathbb{R}^{d \times d}$. Define the data-misfit function as
\[
\Phi(\dv) = \frac{1}{2}\left\|\bfSigma_{\bfeps}^{-\frac{1}{2}}\left(\fhigh(\dv) - \obs\right)\right\|^2\,,
\]
with the norm $\|\cdot\|$. The likelihood function $L: \Zcal \to \mathbb{R}$ is then proportional to
\begin{equation}
L(\obs|\dv) \propto \exp(-\Phi(\dv))\,.
\label{eq:MM:InvProb:LikelihoodHi}
\end{equation}
An evaluation of the likelihood $L$ entails an evaluation of the high-fidelity model $\fhigh$. The posterior probability density is
\[
p(\dv|\obs) \propto L(\obs|\dv)p_0(\dv)\,.
\]

We note that there are many challenging and important questions regarding the formulation of an inverse problem in the Bayesian setting. For example, the selection of the prior is a delicate and problem-dependent question, see for example \cite[Section~3.3]{KaipoBook} for a discussion of prior models. We do not address those issues here, but note that our stated goal of a multifidelity formulation is to recover a solution of the outer-loop problem (here the inverse problem) that retains the accuracy of the high-fidelity formulation. Thus, the multifidelity approaches described below will inherit the choices made in the high-fidelity formulation of the Bayesian inverse problem.

\begin{algorithm}[t]
\caption{Metropolis-Hastings}\label{alg:MCMC-MH}
\begin{algorithmic}[1]
\Procedure{MetropolisHastings}{$L, p_0, q, \nrpts$}
\State Choose a starting point $\dv_0$
\For{$i = 1, \dots, \nrpts$}
\State Draw candidate $\dv^*$ from proposal $q(\cdot | \dv_{i - 1})$
\State Compute acceptance probability
\[
\alpha(\dv_{i - 1}, \dv^*) = \min\left\{1, \frac{q(\dv_{i - 1} | \dv^*)L(\obs|\dv^*)p_0(\dv^*)}{q(\dv^*|\dv_{i - 1})L(\obs|\dv_{i - 1})p_0(\dv_{i-1})}\right\}
\]
\State Set the sample $\dv_i$ to
\[
\dv_i = \begin{cases}
\dv^*\,,\qquad &\text{with probability }\alpha(\dv_{i - 1}, \dv^*)\,,\\
\dv_{i - 1}\,,\qquad &\text{with probability }1 - \alpha(\dv_{i - 1}, \dv^*)
\end{cases}
\]
\EndFor
\State \Return $\dv_1, \dots, \dv_{\nrpts}$
\EndProcedure
\end{algorithmic}
\end{algorithm}

\paragraph{Exploring the posterior}
The solution of the inference problem is explored by drawing samples from the posterior distribution. The posterior samples can then be used to estimate the input with the maximum posterior density and to estimate expectations of functions of interest $h: \Zcal \to \mathbb{R}$ with respect to the posterior distribution. For example, one could be interested in the expected value of $h$ over the posterior distribution
\begin{equation}
\mathbb{E}[h] = \int_{\Zcal} h(\dv)p(\dv|\obs)\mathrm d\dv\,.
\label{eq:MM:InvProb:ExpH}
\end{equation}
MCMC methods are a popular way to sample from the posterior distribution, which have been studied extensively \cite{roberts1997,haario_adaptive_2001,Marzouk20091862,GirolamiHamiltonianMC,cotter2013,GhattasStochNewton,0266-5611-30-11-114015,Cui2016109}; see also the books \cite{GilksMCMC,MCBookLiu}. One example of an MCMC method is the Metropolis-Hastings algorithm, which is an iterative scheme that draws candidate samples from a proposal distribution and then accepts the candidate sample with a probability that depends on the ratio of the posterior at the current candidate sample and the posterior at the sample of the previous iteration.

\paragraph{Metropolis-Hastings algorithm}
Algorithm~\ref{alg:MCMC-MH} summarizes the Metropolis-Hastings algorithm. Inputs are the likelihood $L$, the prior density $p_0$, a proposal density $q$, and the number of samples $\nrpts$. The proposal density $q$ is used to draw the next candidate sample. A typical choice for the proposal density is a Gaussian distribution that is centered at the sample of the previous iteration. In each iteration $i = 1, \dots, \nrpts$, a candidate sample $\dv^*$ is drawn from the proposal $q(\cdot | \dv_{i - 1})$ that depends on the sample $\dv_{i - 1}$ of the previous iteration. The candidate sample $\dv^*$ is accepted $\dv_i = \dv^*$ with the probability $\alpha(\dv^*, \dv_{i - 1})$, which depends on the ratio of the likelihood at the candidate sample $\dv^*$ and the sample $\dv_{i - 1}$ of the previous iteration. If the candidate sample is rejected, then $\dv_i = \dv_{i - 1}$. Algorithm~\ref{alg:MCMC-MH} returns the samples $\dv_1, \dots, \dv_{\nrpts}$.

\paragraph{Metropolis-Hastings algorithm in practice}
Several techniques are necessary to make the Metropolis-Hastings algorithm practical. For example, the samples generated in the first iterations are usually discarded (burn-in) because they are strongly influenced by the initial starting point $\dv_0$. We refer to the literature \cite{RobertBook} for more details and further practical considerations.

\paragraph{Efficiency of MCMC sampling}
Once samples $\dv_1, \dots, \dv_{\nrpts}$ are drawn with an MCMC algorithm, they can be used to, e.g., estimate the expectation \eqref{eq:MM:InvProb:ExpH} as
\[
\bar{h} = \frac{1}{\nrpts}\sum_{i = 1}^{\nrpts} h(\dv_i)\,.
\]
MCMC generates correlated samples $\dv_1, \dots, \dv_{\nrpts}$. The efficiency of MCMC sampling can therefore be measured by the effective sample size for a given computational budget with respect to the estimator $\bar{h}$. To define the effective sample size, consider the samples $\dv_1, \dv_2, \dv_3, \dots$ drawn with an MCMC algorithm and define the integrated autocorrelation time as
\[
\tau_{\text{int}}(h) = \frac{1}{2} + \sum_{j = 1}^{\infty}\rho_j\,,
\]
where $\rho_j = \operatorname{corr}(h(\dv_1), h(\dv_{j + 1}))$ is the correlation between $h(\dv_1)$ and $h(\dv_{j + 1})$. The effective sample size $\nrpts_{\text{eff}}(h)$ is
\[
\nrpts_{\text{eff}}(h) = \frac{\nrpts}{2 \tau_{\text{int}}(h)}\,,
\]
such that
\[
\Var[\bar{h}] \approx \frac{\Var[h]}{\nrpts_{\text{eff}}(h)}\,,
\]
see \cite[p.~125f]{MCBookLiu} for a detailed derivation and further references.
This means that there are at least two ways to improve the efficiency of sampling with MCMC \cite{cui_data-driven_2014}: (1)~Increase the effective sample size for a given number of MCMC iterations with, e.g., adaptive MCMC \cite{gilks_adaptive_1994,gilks_adaptive_1998,roberts_coupling_2007,haario_adaptive_2001}. (2)~Increase the number of MCMC iterations for a given computational budget, so that more samples can be generated for a given budget with, e.g., two-stage MCMC \cite{christen_markov_2005,Fox97samplingconductivity}.

\subsection{Two-stage Markov chain Monte Carlo}
\label{sec:MM:InvProb:TwoStageMCMC}
Two-stage MCMC methods aim to increase the number of MCMC iterations for a given computational budget. In many applications, the Metropolis-Hastings algorithm, and MCMC in general, requires many iterations to produce an acceptable effective sample size. Each iteration means a likelihood evaluation, which means a high-fidelity model $\fhigh$ evaluation in the case of the likelihood $L$ as defined in \eqref{eq:MM:InvProb:LikelihoodHi}. Two-stage MCMC methods employ delayed acceptance or delayed rejection strategies that use multifidelity models to reduce the number of samples evaluated using the expensive high-fidelity model.

The work \cite{christen_markov_2005,Fox97samplingconductivity} proposes a two-stage delayed acceptance MCMC sampling. The candidate sample $\dv^*$ has to be accepted with the likelihood induced by a low-fidelity model first, before $\dv^*$ is passed on to be either accepted or rejected with the likelihood induced by the high-fidelity model.

Multi-stage MCMC methods based on delayed rejection, in contrast to delayed acceptance, have been proposed in \cite{tierney_adaptive_1999,mira_metropolis-hastings_2001,green_delayed_2001}. A candidate sample that is rejected in one stage is retried with a different proposal distribution at a subsequent stage. Ref.~\cite{tierney_adaptive_1999} suggests using an independence sampler based on a density that is thought to be a good approximation of the posterior density on the first stage. In case the density corresponding to the independence sampler is a poor approximation, a random walk is used on the second stage. Another approach discussed in \cite{tierney_adaptive_1999} is to use a proposal based on a local quadratic approximation of the log posterior density, which is restricted to smaller and smaller neighborhoods of the current point in a trust-region fashion. In \cite{green_delayed_2001}, the proposal at the first stage is a normal distribution. The proposal at the second stage is a normal distribution with the same mean but with a higher variance. In principle, the proposals in the context of delayed rejection MCMC, just as in delayed acceptance MCMC, can also be derived from models with different fidelities and costs (although we are not aware of such implementations in the literature that pre-date the work by \cite{christen_markov_2005,Fox97samplingconductivity}).

\paragraph{Two-stage delayed acceptance MCMC algorithm}
The two-stage MCMC method introduced in \cite{christen_markov_2005} is summarized in Algorithm~\ref{alg:MCMC-2Stage}. Inputs are the likelihood $L^{\hi}$ and $L^{\lo}$, corresponding to the high- and low-fidelity model respectively, the prior density $p_0$, the proposal $q$, and the number of samples $\nrpts$. Consider one of the iterations $i = 1, \dots, \nrpts$. The first stage of Algorithm~\ref{alg:MCMC-2Stage} (lines 4--6) proceeds as the Metropolis-Hastings algorithm with a candidate sample $\dv^{\lo}$ drawn from the proposal distribution, except that the likelihood of the low-fidelity model is used, instead of the likelihood of the high-fidelity model. The result of the first stage is a sample $\dv^{\hi}$, which either is $\dv^{\hi} = \dv^{\lo}$ (accept) or $\dv^{\hi} = \dv_{i - 1}$ (reject). In the second stage of the algorithm (lines 7--9), the high-fidelity model is used to either accept or reject the candidate sample $\dv^{\hi}$ of the first stage. If the sample $\dv^{\hi}$ is accepted in the second stage, then the algorithm sets $\dv_i = \dv^{\hi}$. If the sample $\dv^{\hi}$ is rejected in the second stage, then $\dv_i = \dv_{i - 1}$. Note that in case the first stage rejected the sample $\dv^{\lo}$, i.e., $\dv^{\hi} = \dv_{i - 1}$, no high-fidelity model evaluation is necessary in the second stage because the high-fidelity model output at $\dv_{i - 1}$ is available from the previous iteration.
The proposal distribution $Q$ in the second stage depends on the likelihood $L^{\lo}$ of the low-fidelity model. Note that $\delta_{\dv_{i - 1}}$ is the Dirac mass at $\dv_{i - 1}$. Algorithm~\ref{alg:MCMC-2Stage} returns the samples $\dv_1, \dots, \dv_{\nrpts}$.

\paragraph{Efficiency of two-stage delayed acceptance MCMC}
The key to the efficiency of Algorithm~\ref{alg:MCMC-2Stage} is that no high-fidelity model evaluation is necessary in the second stage if the candidate sample $\dv^{\lo}$ was rejected at the first stage. Since the rejection rate of MCMC is typically high, many high-fidelity model evaluations are saved by the two-stage MCMC. We refer to \cite{christen_markov_2005} for an asymptotic analysis and the convergence properties of the two-stage MCMC. The two-stage MCMC uses a model management based on filtering because only candidate samples accepted with the low-fidelity model are passed on to the high-fidelity model, see Section~\ref{sec:MM:strategies}.

\paragraph{Other multifidelity MCMC approaches}
In \cite{efendiev_preconditioning_2006}, two-stage MCMC is seen as a preconditioned Metropolis-Hastings algorithm. The low-fidelity model in \cite{efendiev_preconditioning_2006} is a coarse-scale model of a high-fidelity multiscale finite volume model. We also mention \cite{higdon_bayesian_2002}, where multiple MCMC chains from coarse- (low-fidelity) and fine-scale (high-fidelity) models are coupled using a product chain. The work \cite{MultiMCMC} couples multilevel Monte Carlo and MCMC to accelerate the estimation of expected values with respect to a posterior distribution. The low-fidelity models form a hierarchy of coarse-grid approximations of the high-fidelity model.

\begin{algorithm}[t]
\caption{Two-stage MCMC}\label{alg:MCMC-2Stage}
\begin{algorithmic}[1]
\Procedure{2StageMCMC}{$L^{\hi}, L^{\lo}, p_0, q, \nrpts$}
\State Choose a starting point $\dv_0$
\For{$i = 1, \dots, \nrpts$}
\State Draw candidate $\dv^{\lo}$ from proposal $q(\cdot | \dv_{i - 1})$
\State Compute acceptance probability
\[
\alpha^{\lo}(\dv_{i - 1}, \dv^{\lo}) = \min\left\{1, \frac{q(\dv_{i - 1} | \dv^{\lo})L^{\lo}(\obs|\dv^{\lo})p_0(\dv^{\lo})}{q(\dv^{\lo}|\dv_{i - 1})L^{\lo}(\obs|\dv_{i - 1})p_0(\dv_{i-1})}\right\}
\]
\State Set the candidate sample $\dv^{\hi}$ to
\[
\dv^{\hi} = \begin{cases}
\dv^{\lo}\,,\qquad &\text{with probability }\alpha^{\lo}(\dv_{i - 1}, \dv^{\lo})\,,\\
\dv_{i - 1}\,,\qquad &\text{with probability }1 - \alpha^{\lo}(\dv_{i - 1}, \dv^{\lo})
\end{cases}
\]
\State Set the distribution $Q$ to
\begin{multline*}
Q(\dv | \dv_{i - 1}) = \alpha^{\lo}(\dv_{i - 1}, \dv)q(\dv|\dv_{i - 1}) + \\\left(1 - \int_{\Zcal}\alpha^{\lo}(\dv_{i - 1}, \dv)q(\dv | \dv_{i - 1})\mathrm d\dv\right)\delta_{\dv_{i - 1}}(\dv)
\end{multline*}
\State Compute acceptance probability
\[
\alpha^{\hi}(\dv_{i - 1}, \dv^{\hi}) = \min\left\{1, \frac{Q(\dv_{i - 1} | \dv^{\hi})L^{\hi}(\obs|\dv^{\hi})p_0(\dv^{\hi})}{Q(\dv^{\hi}|\dv_{i - 1})L^{\hi}(\obs|\dv_{i - 1})p_0(\dv_{i-1})}\right\}
\]
\State Set sample $\dv_i$ to
\[
\dv_i = \begin{cases}
\dv^{\hi}\,,\qquad &\text{with probability }\alpha^{\hi}(\dv_{i - 1}, \dv^{\hi})\,,\\
\dv_{i - 1}\,,\qquad &\text{with probability }1 - \alpha^{\hi}(\dv_{i - 1}, \dv^{\hi})
\end{cases}
\]
\EndFor
\State \Return $\dv_1, \dots, \dv_{\nrpts}$
\EndProcedure
\end{algorithmic}
\end{algorithm}

\subsection{Markov chain Monte Carlo with adaptive low-fidelity models}
\label{sec:MM:InvProb:Adaptive}
In \cite{cui_data-driven_2014} an algorithm for combining high- and low-fidelity models in MCMC sampling is presented. The low-fidelity model is used in the first step of a two-stage MCMC approach to increase the acceptance rate of candidates in the second step, where the high-fidelity model is used to either accept or reject the sample. Additionally, the high-fidelity model outputs computed in the second step are used to adapt the low-fidelity model. In that sense, the approach uses a model management based on a combination of adaptation and filtering.

The low-fidelity model is a projection-based reduced model in \cite{cui_data-driven_2014}, see Section~\ref{sec:intro:models}. The low-fidelity model is constructed in an offline phase from an initial reduced basis. At each MCMC iteration, the low-fidelity model is used in a first stage to generate a certain number of samples with the Metropolis-Hastings algorithm. The goal is to generate so many samples with the low-fidelity model that the initial sample and the last sample are almost uncorrelated. Then, in the second stage of the MCMC iteration, the last sample generated with the low-fidelity model is used as a candidate sample and the acceptance probability is computed using the high-fidelity model. Similarly to the two-stage MCMC methods discussed above, this algorithm aims to increase the acceptance probability at the second stage; however, the high-fidelity model output is used to improve the low-fidelity model after a sample has been accepted, and is not discarded. In that way, the authors of \cite{cui_data-driven_2014} improve the low-fidelity model during the MCMC iterations, and consequently increase the acceptance probability in the second stage of the MCMC iterations. The authors show that under certain technical conditions, the Hellinger distance $d_{\text{Hell}}(p_{\hi}, p_{\lo})$ between the posterior distribution $p_{\hi}$ corresponding to the high-fidelity model $\fhigh$ and the approximate posterior distribution $p_{\lo}$ induced by the reduced model $\fl$ decreases as the reduced model $\fl$ is improved during the MCMC iterations. Thus, for any $\epsilon > 0$, the Hellinger distance between the high-fidelity and low-fidelity posterior distribution can be bounded
\[
d_{\text{Hell}}(p_{\hi}, p_{\lo}) \leq \epsilon
\]
after a sufficiently large number of MCMC iterations. The resulting low-fidelity model is data-driven because it uses information provided by the observation $\obs$, rather than only prior information. Using samples from this adaptive two-stage MCMC approach yields unbiased Monte Carlo estimators, see the analysis in \cite{cui_data-driven_2014}. 

\subsection{Bayesian estimation of low-fidelity model error}
\label{sec:MM:InvProb:Error}
In \cite{KaipoBook,kaipio_statistical_2007}, a multifidelity approach for Bayesian inference is presented that relies on a Bayesian approximate error model of the low-fidelity model derived from the high-fidelity model. The inference is performed with the low-fidelity model but an additional term is introduced that quantifies the error between the low-fidelity and the high-fidelity model. The error term relies on the high-fidelity model and is used to correct for the error introduced by the low-fidelity model. The approach follows the Bayesian paradigm and considers the error term as a random variable, which becomes another noise term in the Bayesian formulation of the inverse problem. Within our categorization, this approach uses a model management based on model adaptation.

In \cite{KaipoBook,kaipio_statistical_2007}, the relationship $\obs = \fhigh(\dv) + \bfeps$ is formulated as
\[
\obs = \fl(\dv) + \left(\fhigh(\dv) - \fl(\dv)\right) + \bfeps\,,
\]
which can be rewritten as
\[
\obs = \fl(\dv) + \bfe(\dv) + \bfeps\,,
\]
with the error term $\bfe(\dv)$. If $\fhigh$ and $\fl$ are linear in $\dv$, the sum of the covariance of the error $\bfe(\dv)$ and noise term $\bfeps$ can derived explicitly under certain assumptions \cite{KaipoBook,kaipio_statistical_2007}. If $\fhigh$ and $\fl$ are nonlinear, a Gaussian approximation of the covariance is proposed in \cite{KaipioErrorDOT}.

The Bayesian approximate error model is used in \cite{KaipioErrorDOT,kolehmainen_approximation_2012} for a diffuse optical tomography problem, where the low-fidelity model corresponds a coarse-grid finite-element approximation and the high-fidelity model is a fine-grid approximation. In \cite{tarvainen_approximation_2010,tarvainen_corrections_2010}, diffuse optical tomography is considered with a simplified-physics low-fidelity model. The error introduced by the low-fidelity model is corrected with the Bayesian approximate error modeling approach. Correcting for errors due to truncation of the computational domain and unknown domain shape are investigated in \cite{lehikoinen_approximation_2007,nissinen_compensation_2011}. The Bayesian approximate error model approach was verified on real data in \cite{nissinen_bayesian_2008}. Extensions to time-dependent problems are presented in \cite{huttunen_approximation_2007,huttunen_approximationB_2007}.

We note that estimating the error term $\bfe$ is related to estimating the model inadequacy of the low-fidelity model with respect to the high-fidelity model. Often, Gaussian-process models are used to attempt to correct model inadequacy \cite{kennedy_bayesian_2001,kennedy_bayesian_2002,higdon_combining_2004,oakley_bayesian_2002,bayarri_framework_2007}. The corrected low-fidelity model is typically cheap to evaluate, because the Gaussian-process model is cheap, and so the corrected low-fidelity model can be used in MCMC sampling with the aim of decreasing the costs per MCMC iteration compared to using the high-fidelity model. The work \cite{KaipoBook,kaipio_statistical_2007} can be seen as dealing with model inadequacy in the specific context of inverse problems.

%% file: opti.tex
\section{Multifidelity model management in optimization}
\label{sec:MM:Optimization}
The goal of optimization is to find an input that leads to an optimal model output with respect to a given objective function. Optimization is typically formulated as an iterative process that requires evaluating a model many times. Using only an expensive high-fidelity model is often computationally prohibitive. We discuss multifidelity methods that leverage low-fidelity models for speeding up the optimization while still resulting in a solution that satisfies the optimality conditions associated with the high-fidelity model. Section~\ref{sec:MM:Optimization:Opti} formalizes the optimization task. Section~\ref{sec:MM:Optimization:Global} discusses global multifidelity optimization, which searches over the entire feasible domain, and Section~\ref{sec:MM:Optimization:Local} reviews local multifidelity optimization, which searches for a locally optimal solution, typically in a neighborhood of an initial input.

\subsection{Optimization using a single high-fidelity model}
\label{sec:MM:Optimization:Opti}
We consider the setting of unconstrained optimization. Given is the high-fidelity model $\fhigh: \Zcal \to \Ycal$, with the $d$-dimensional input $\dv \in \Zcal$. The goal is to find an input $\dv^* \in \Zcal$ that solves
\begin{equation}
\min_{\dv \in \Zcal} \fhigh(\dv)\,.
\label{eq:MM:opt:problem}
\end{equation}
Typically, the optimal input $\dv^*$ is obtained in an iterative way, where a sequence of inputs $\dv_1, \dv_2, \dots$ is constructed such that this sequence $(\dv_i)$ converges to the optimal input $\dv^*$ in a certain sense \cite{Wright}.

\paragraph{Local and global optimization}
We distinguish between local and global optimization approaches. Global optimization searches over the entire feasible domain $\Zcal$ for a minimizer $\dv^*$, whereas local optimization terminates when a local optimum is found. Local optimization thus searches in a neighborhood of an initial point $\dv \in \Zcal$. Global methods typically do not require the gradient of $\fhigh$, which may be advantageous in situations where the model is given as a black box and where approximations of the gradient contain significant noise.
The use of gradient information in local methods leads to a more efficient search process that typically uses fewer model evaluations and that is more scalable to problems with high-dimensional input. There are also local methods that do not require gradient information.

We note that the multifidelity solution of constrained optimization problems can be achieved using penalty formulations that convert the constrained problem into an unconstrained one, although some multifidelity methods admit more sophisticated ways to handle approximations of the constraints. Local methods can more easily deal with constraints, whereas global methods often fall back to heuristics (see, e.g., the discussion in \cite{MarchPhD}).

\subsection{Global multifidelity optimization}
\label{sec:MM:Optimization:Global}
A large class of global multifidelity optimization methods uses adaptation as a model management strategy. These methods
search for a minimizer with respect to an adaptively refined low-fidelity model $\fl$. They guarantee that the resulting optimal solution is also a minimizer of the high-fidelity model $\fhigh$ by the way in which the low-fidelity model is adapted throughout the optimization search, using information from high-fidelity model evaluations. A critical task in this class of global multifidelity optimization methods therefore is to balance between exploitation (i.e., minimizing the current low-fidelity model) and exploration (i.e., adapting the low-fidelity model with information from the high-fidelity model).

\subsubsection{Efficient global optimization (EGO)}
EGO with expected improvement is frequently used to balance exploitation and exploration in cases where the low-fidelity model is a kriging model, see \cite{jones_efficient_1998} and Section~\ref{sec:intro:models}. Let $\dv$ be the input that minimizes the low-fidelity model $\fl$ at the current iteration. Sampling at a new point $\dv^{\prime} \in \Zcal$ means evaluating the high-fidelity model at $\dv^{\prime}$ and then adapting the low-fidelity model with information obtained from $\fhigh(\dv^{\prime})$.

EGO provides a formulation for choosing the new sample point $\dv^{\prime}$. It does this by considering the value $\fhigh(\dv^{\prime})$ of the high-fidelity model at the new sampling point $\dv^{\prime}$ to be uncertain before the high-fidelity model is evaluated at $\dv^{\prime}$. The uncertainty in $\fhigh(\dv^{\prime})$ is modeled as a Gaussian random variable $Y$ with mean and standard deviation given by the kriging (low-fidelity) model, see Section~\ref{sec:intro:models}. The improvement at point $\dv^{\prime}$ is then given by $I(\dv^{\prime}) = \max\{\fhigh(\dv) - Y, 0\}$, and the expected improvement is
\[
\mathbb{E}[I(\dv^{\prime})] = \mathbb{E}[\max\{\fhigh(\dv) - Y, 0\}]\,.
\]
The expected improvement at a point $\dv^{\prime}$ can be efficiently computed by exploiting the cheap computation of the MSE (and therefore standard deviation) of kriging models \cite{ForresterKriging,jones_efficient_1998}. The optimization process is then to start with an initial kriging model, find a new input that maximizes the expected improvement, evaluate the high-fidelity model at this new input, update the kriging model, and iterate. The optimization loop is typically stopped when the expected improvement is less than some small positive number \cite{jones_taxonomy_2001}.

EGO can be made globally convergent  and does not require high-fidelity model derivatives \cite{jones_taxonomy_2001}. However, as discussed in \cite{jones_taxonomy_2001,MarchPhD}, EGO is sensitive to the initial points used for building the kriging model, which means that first a fairly exhaustive search around the initial points might be necessary before a more global search begins. The reason for this behavior is that EGO proceeds in two steps. In the first step, a kriging model is learned, and in the second step, the kriging model is used as if it were correct to determine the next point. We refer to \cite{jones_taxonomy_2001,ForresterKriging} for details and possible improvements. The work in \cite{RajnarayanPhD}, builds on EGO and develops a concept of expected improvement for multiobjective optimization that considers improvement and risk when selecting a new point. We also mention \cite{chaudhuri_efficient_2014}, where EGO is equipped with an adaptive target method that learns from previous iterations how much improvement to expect in the next iteration.

\subsubsection{Other approaches to combine multifidelity models in the context of global optimization}
In the work \cite{allaire_mathematical_2014}, information from multiple kriging models is fused. Each model in the multifidelity hierarchy is approximated by a kriging model, and the random variables representing the kriging models are fused following the technique introduced in \cite{winkler_combining_1981}. 
The authors of \cite{Goel2007} propose to use a weighted average of an ensemble of low-fidelity models in an optimization framework. The weights are derived from the errors of the low-fidelity models.

In \cite{gramacy_adaptive_2009}, Gaussian-process low-fidelity models are adapted via active learning, instead of EGO. The approach samples the input space and efficiently fits local Gaussian-process models. These local models are obtained via the treed Gaussian-process implementation \cite{gramacy_bayesian_2008} that partitions the input space into several regions and builds separate Gaussian-process models in each region. The approach in \cite{gramacy_adaptive_2009} is particularly suited to the case where sampling the high-fidelity model entails large-scale, parallel computations. In \cite{gramacy_modeling_2016}, Gaussian-process models and expected improvement are developed for constrained optimization. The augmented Lagrangian is used to reduce the constrained optimization problem into an unconstrained one, for which Gaussian-process models can be derived.

There are multifidelity optimization methods based on pattern search \cite{Hooke:1961:DSS:321062.321069,PSConvergence} that use filtering as model management strategy. For example, the multifidelity pattern search algorithm presented in \cite{booker_rigorous_1999} uses low-fidelity models to provide additional search directions, while preserving the convergence to a minimizer of the high-fidelity model. In \cite{taddy_bayesian_2009}, a Bayesian approach to pattern search with Gaussian-process models is presented. The low-fidelity model is used to guide the pattern search. We also refer to \cite{Queipo20051} for a discussion of pattern search algorithms with low-fidelity models.

\subsection{Local multifidelity optimization}
\label{sec:MM:Optimization:Local}
Local multifidelity optimization methods in the literature typically use adaptation as a model management strategy. One class of approaches uses direct adaptation of a low-fidelity model, using high-fidelity evaluations that are computed as the optimization proceeds. In each optimization iteration, the high-fidelity model is evaluated at the minimizer of the low-fidelity model. The corresponding high-fidelity model output is then used to adapt the low-fidelity model. Convergence to the minimizer of the high-fidelity model can be guaranteed in limited situations, see \cite{Barthelemy1993,Queipo20051,ForresterKriging} for details. Another class of approaches performs optimization on a corrected low-fidelity model, where the corrections are computed from high-fidelity model evaluations as the optimization proceeds. These latter approaches typically use a trust-region framework to manage the corrections, which we discuss in more detail now.

\subsubsection{Multifidelity trust-region methods}
A classical way of exploiting a low-fidelity model in an optimization framework is to optimize over the low-fidelity model in a trust region---that is, to solve an optimization subproblem using the low-fidelity model, but to restrict the size of the step to lie within the trust region. A classical example is to derive a quadratic approximation of the high-fidelity model with its gradient and Hessian at the center of the trust region. The size of the trust region is then determined depending on the approximation quality of the quadratic approximation. We refer to \cite{conn_trust_2000} for an introduction, theory, and implementation of these classical trust region methods. Early formulations of trust region ideas can be found in the papers by Levenberg \cite{levenberg_method_1944} and Marquardt \cite{marquardt_algorithm_1963}, with the papers by Powell \cite{PowellConv3,PowellConv2,PowellConv} providing the important contribution of establishing and proving convergence of trust-region algorithms.

Classical trust-region methods use quadratic approximations of the objective in the trust region. The work \cite{alexandrov_trust-region_1998} established a multifidelity trust region framework for more general low-fidelity models. In particular, \cite{alexandrov_trust-region_1998} formulates a first-order consistency requirement on the low-fidelity model. The first-order consistency requirement is that the low- and the high-fidelity model have equal value and gradient at the center of the trust region. This consistency requirement ensures that the resulting multifidelity optimization algorithm converges to an optimal solution of the high-fidelity model.

\begin{algorithm}[t]
\caption{A multifidelity trust-region algorithm}\label{alg:Opti-TR}
\begin{algorithmic}[1]
\Procedure{TrustRegion}{$\fhigh, \eta_1, \eta_2, \gamma_1, \gamma_2, \Delta^*$}
\State Choose a starting point $\dv_1$ and initial step size $\Delta_1 > 0$
\For{$i = 1, 2, 3, \dots$ until convergence}
\State Construct low-fidelity model $\fl^{(i)}$ with
\[
\hspace*{2cm}\fl^{(i)}(\dv_i) = \fhigh(\dv_i)\,,\text{ and } \nabla\fl^{(i)}(\dv_i) = \nabla \fhigh(\dv_i)
\]
\State Find $\boldsymbol s_i \in \Zcal$ that solves
\begin{align*}
\hspace*{2cm}\min \quad & \fl^{(i)}(\dv_i + \boldsymbol s_i)\\
\hspace*{2cm}\text{s.t.}\quad & \|\boldsymbol s_i\| \leq \Delta_i
\end{align*}
\State Update point
\[
\hspace*{2cm}\dv_{i + 1} = \begin{cases}
\dv_i + \boldsymbol s_i\,,\quad \text{ if } \fhigh(\dv_i + \boldsymbol s_i) < \fhigh(\dv_i)\,,\\
\dv_i\,,\quad \text{otherwise}
\end{cases}
\]
\State Compare actual and estimated decrease in $\fhigh$
\[
\hspace*{2cm}\gamma = \frac{\fhigh(\dv_i) - \fhigh(\dv_i + \boldsymbol s_i)}{\fhigh(\dv_i) - \fl^{(i)}(\dv_i + \boldsymbol s_i)}
\]
\State Update step size
\[
\hspace*{2cm}\Delta_{i + 1} = \begin{cases}
\eta_1\|\boldsymbol s_i\|\,,\qquad & \gamma < \gamma_1\,,\\
\min\{\eta_2\Delta_i,\Delta^*\}\,,\qquad &\gamma > \gamma_2\,,\\
\|\boldsymbol s_i\|\,,\qquad &\text{otherwise}
\end{cases}
\]
\EndFor
\State \Return $\dv_1, \dv_2, \dv_3, \dots$
\EndProcedure
\end{algorithmic}
\end{algorithm}

\subsubsection{Multifidelity trust-region algorithm}
The multifidelity trust-region approach of \cite{alexandrov_trust-region_1998} is summarized in Algorithm~\ref{alg:Opti-TR}. Inputs are the high-fidelity model $\fhigh$, the parameters $\eta_1, \gamma_1, \eta_2, \gamma_2 > 0$ that control the expansion and contraction of the trust region, and an upper bound $\Delta^* > 0$ on the step size. In each iteration of the loop, an update to the current point is proposed and the trust-region is either expanded or contracted. Note that we have omitted stopping criteria in Algorithm~\ref{alg:Opti-TR} for the ease of exposition. In each iteration $i = 1, 2, 3, \dots$, a low-fidelity model $\fl^{(i)}$ is constructed so that it satisfies the first-order consistency requirement. Thus, the low-fidelity model output $\fl^{(i)}(\dv_i) = \fhigh(\dv_i)$ equals the high-fidelity model output at the current point $\dv_i$, and the gradient of the low-fidelity model $\nabla\fl^{(i)}(\dv_i) = \nabla\fhigh(\dv_i)$ equals the gradient of the high-fidelity model. It is shown in \cite{alexandrov_trust-region_1998} that with an additive or multiplicative correction an arbitrary low-fidelity model can be adjusted to satisfy the first-order consistency requirement. The first-order consistency requirement can also be established with the scaling approach introduced in \cite{haftka_combining_1991}. Then, the step $\boldsymbol s_i$ is computed and accepted if the point $\dv_i + \boldsymbol s_i$ decreases the high-fidelity model $\fhigh$ with respect to the current point $\dv_i$, i.e., if $\fhigh(\dv_i + \boldsymbol s_i) < \fhigh(\dv_i)$. The size of the trust region is updated depending on the ratio of the actual decrease $\fhigh(\dv_i) - \fhigh(\dv_i + \boldsymbol s_i)$ to the estimated decrease $\fhigh(\dv_i) - \fl^{(i)}(\dv_i + \boldsymbol s_i)$. The parameters $\gamma_1$ and $\gamma_2$ control when to shrink and when to expand the trust region, respectively. The parameters $\eta_1$ and $\eta_2$ define the size of the contracted and expanded trust region, respectively. The algorithm returns the points $\dv_1, \dv_2, \dv_3, \dots$.

\subsubsection{Generalized multifidelity trust-region methods}
In Algorithm~\ref{alg:Opti-TR}, the low-fidelity model is corrected in each iteration. In \cite{Arian2002,Fahl2003}, the trust-region POD (TRPOD) is proposed, which uses a projection-based POD reduced model in conjunction with a trust-region optimization framework and adapts the POD reduced model in each optimization iteration. Similarly, in \cite{Bergmann20087813,zahr_progressive_2015}, projection-based models are adapted within similar trust region frameworks. The work \cite{yue2013} uses error bounds for interpolatory reduced models to define the trust region and refinement of the reduced models to guarantee convergence to the optimality conditions associated with the high-fidelity model. In \cite{fischer_bayesian_2017}, a correction model is used to correct low-fidelity model outputs in a trust-region model management optimization scheme. The correction model is constructed in the trust region and utilized to calculate additive and multiplicative adjustments to the low-fidelity model. The work \cite{robinson_surrogate-based_2008} develops a trust-region multifidelity approach for combining low- and high-fidelity models that are defined over different design spaces. For example, the number of design variables used with the low-fidelity model can be different than the number of design variables with the high-fidelity model. Space mappings are derived that map between the design spaces. The space mappings are corrected to obtain a provably convergent method.

If gradients are unavailable, or too costly to approximate, then the framework developed by \cite{alexandrov_trust-region_1998} relying on the first-order consistency of the low-fidelity model cannot be applied. In \cite{march_constrained_2012}, a gradient-free multifidelity trust-region framework with radial basis function interpolants as low-fidelity models is introduced, building on the gradient-free trust-region methods of \cite{Conn2009,ConnDFree} and the tailored radial basis function modeling of \cite{Wild2}. This leads to an error interpolation that makes the low-fidelity model satisfy the sufficient condition to prove convergence to a minimizer of the high-fidelity model.

\subsection{Multifidelity optimization under uncertainty}
In their simplest forms, optimization under uncertainty formulations consider a similar problem to (\ref{eq:MM:opt:problem}), but where now the objective function $\fhigh$ incorporates one or more statistics that in turn depend on underlying uncertainty in the problem. For example, a common objective function is a weighted sum of expected performance and performance standard deviation. Thus, each optimization iteration embeds an uncertainty quantification loop (e.g., Monte Carlo sampling or stochastic collocation) over the uncertain parameters. \cite{ng_multifidelity_2014}  uses the control variate-based multifidelity Monte Carlo method of Section~\ref{sec:MM:UQ:ContVar} to derive a multifidelity optimization under uncertainty method that provides substantial reductions in computational cost using a variety of low-fidelity model types. In  \cite{ng_monte_2015}, it is shown that evaluations computed at previous optimization iterations form an effective and readily-available low-fidelity model that can be exploited by this control-variate formulation in the optimization under uncertainty setting. The work in \cite{Shah2015} uses a combination of polynomial chaos stochastic expansions and corrections based on coarse-grid approximations to formulate a multifidelity robust optimization approach.

We highlight optimization under uncertainty as an important target area for future multifidelity methods. It is a computationally demanding process, but one with critical importance to many areas, such as engineering design.

%% file: outlook.tex
\section{Conclusions and outlook}
\label{sec:Outlook}

As can be seen from the vast literature surveyed in this paper, multifidelity methods have begun to have impact across diverse outer-loop applications in computational science and engineering. Multifidelity methods have been used for more than two decades to accelerate solution of optimization problems. Their application in uncertainty quantification is more recent, but appears to offer even more opportunity for computational speedups, due to the heavy computational burden typically associated with uncertainty quantification tasks such as Monte Carlo and MCMC sampling.

This paper highlights the broad range of multifidelity approaches that exist in the literature. These approaches span many types of low-fidelity models as well as many specific strategies for achieving the multifidelity model management. We attempt to bring some perspective on the similarities and differences across methods, as well as their relative advantages, by categorizing methods that share a common theoretical foundation. We discuss methods to create low-fidelity models according to the three areas of simplified models, projection-based models, and data-fit models. We categorize multifidelity model management methods as being based on adaptation, fusion, and filtering. In most settings, one can flexibly choose a combination of model management strategy and low-fidelity model type, although---as is always the case in computational modeling---bringing to bear knowledge of the problem structure helps to make effective decisions on the low-fidelity modeling and multifidelity model management strategies that are best suited to the problem at hand. We note that this paper focused on outer-loop applications and therefore mentioned only briefly multifidelity methods that do not directly exploit the outer-loop setting.

Multifidelity methods have advanced considerably, especially in the past decade. Yet a number of important challenges remain. In almost all existing multifidelity methods (and in the presentation of material in this paper), the assumption is that the high-fidelity model represents some ``truth." This ignores the important fact that the output of the high-fidelity model is itself an approximation of reality. Even in the absence of uncertainties in the inputs, all models, including the high-fidelity model, are inadequate \cite{kennedy_bayesian_2001}. Furthermore, the relationships among different models may be much richer than the linear hierarchy assumed by existing multifidelity methods. One way multifidelity methods can account for model inadequacy is by fusing outputs from multiple models. Model inadequacy is often quantified by probabilities, which describe the belief that a model yields the true output. Techniques for assigning model probabilities reach from expert opinion \cite{zio_two_1996,reinert_including_2006} to statistical information criteria \cite{link_model_2006,BurnhamMultimodel} to quantifying the difference between experimental data and model outputs \cite{park_bayesian_2010}. The probabilities are then used to fuse the model outputs with, e.g., Bayesian model averaging \cite{LeamerBook} or adjustment factors \cite{reinert_including_2006}. The fundamental work by Kennedy and O'Hagan \cite{kennedy_bayesian_2001} treats model inadequacy in a Bayesian setting and led to a series of papers on correcting model inadequacy and validating models \cite{kennedy_bayesian_2002,higdon_combining_2004,oakley_bayesian_2002,bayarri_framework_2007}. The work \cite{bayarri_framework_2007} introduces a six-step process to validate models. In \cite{allaire_mathematical_2014}, Bayesian estimation is employed to fuse model outputs together. Incorporating such approaches into multifidelity model management strategies for outer loop applications remains an important open research challenge.

Another important challenge, already discussed briefly in this paper, is to move beyond methods that focus exclusively on models, so that decision-makers can draw on a broader range of available information sources. Again, some of the foundations exist in the statistical literature, such as \cite{kennedy_predicting_2000} which derives models by fusing multiple information sources such as experiments, expert opinions, lookup tables, and computational models. In drawing on various information sources, multifidelity model management strategies must be expanded to address not just which information source to evaluate when, but also {\em where} (i.e., at what inputs) to evaluate the information source. Relevant foundations to address this challenge include the experimental design literature and value of information analysis \cite{poloczek2016}.